\documentclass[12pt]{article}

\newcommand{\blind}{1}

\usepackage{float}
\usepackage{amsmath}
\usepackage{graphicx, subfigure}
\usepackage[round]{natbib}
\usepackage{geometry}
\usepackage{enumerate}
\usepackage{tikz}
\usepackage[english]{babel}
\usepackage{longtable}
\usepackage{color}
\usepackage{amssymb, amsmath, amsthm}
\usepackage{multirow}
\usepackage{paralist}
\usepackage{authblk}
\usepackage{setspace}
\usepackage{dsfont}
\usepackage{enumerate}
\usepackage[titletoc,title]{appendix}
\usepackage{float}
\usepackage{dsfont}
\usepackage{amsmath}
\usepackage{amsfonts}
\usepackage{multirow}
\usepackage{enumerate}

\usepackage{hyperref}
\hypersetup{
    colorlinks=true,
    linkcolor=blue,
    filecolor=magenta,
    urlcolor=blue,
    citecolor=blue,
}
\providecommand{\keywords}[1]{\textbf{\textit{Keywords:}} #1}

\newtheorem{theorem}{Theorem}
\newtheorem{lemma}{Lemma}

\newcommand{\indep}{\mbox{$\perp\!\!\!\perp$}}

\DeclareMathOperator{\logit}{logit}

\newcommand{\equivReplacedByEquality}{=}

\newcommand{\one}{\mathds{1}}

\newcommand{\sigmakm}{\sigma_{\footnotesize \mbox{km}}}
\newcommand{\sigmaadjipw}{\sigma_{\footnotesize \mbox{adj,ipw}}}
\newcommand{\hatsigmaadjeff}{\hat{\sigma}_{\footnotesize \mbox{adj,eff}}}
\newcommand{\sigmaadjeff}{\sigma_{\footnotesize \mbox{adj,eff}}}

\newcommand{\Dkm}{D_{\footnotesize \mbox{km}}}
\newcommand{\Skmat}{\hat{S}_{\footnotesize \mbox{km}}(t,a)}
\newcommand{\Gkmat}{\hat{G}_{\footnotesize \mbox{km}}(t,a)}
\newcommand{\Skmot}{\hat{S}_{\footnotesize \mbox{km}}(t,1)}
\newcommand{\Skmzt}{\hat{S}_{\footnotesize \mbox{km}}(t,0)}

\newcommand{\thetaTMLE}{\hat{\theta}_{\footnotesize \mbox{adj,eff}}}

\newcommand{\SIPWat}{\hat{S}_{\footnotesize \mbox{ipw}}(t,a)}
\newcommand{\SIPWot}{\hat{S}_{\footnotesize \mbox{ipw}}(t,1)}
\newcommand{\SIPWzt}{\hat{S}_{\footnotesize \mbox{ipw}}(t,0)}

\newcommand{\SADJIPWat}{\hat{S}_{\footnotesize \mbox{adj,ipw}}(t,a)}
\newcommand{\SADJIPWot}{\hat{S}_{\footnotesize \mbox{adj,ipw}}(t,1)}
\newcommand{\SADJIPWzt}{\hat{S}_{\footnotesize \mbox{adj,ipw}}(t,0)}

\newcommand{\thetaIPW}{\hat{\theta}_{\footnotesize \mbox{ipw}}}
\newcommand{\thetaAIPW}{\hat{\theta}_{\footnotesize \mbox{aipw}}}
\newcommand{\thetaADJIPW}{\hat{\theta}_{\footnotesize \mbox{adj,ipw}}}
\newcommand{\thetakm}{\hat{\theta}_{\footnotesize \mbox{km}}}

\title{Improved Precision in the Analysis of Randomized Trials with Survival Outcomes, without Assuming Proportional Hazards}

\if1\blind
\author[1] {Iv\'an  D\'iaz}
\author[2]{Elizabeth Colantuoni}
\author[3]{Daniel F. Hanley}
\author[2]{Michael Rosenblum}
\affil[1]{\small Division of Biostatistics and Epidemiology, Weill Cornell Medicine, New York, NY, USA}
\affil[2]{\small Department of Biostatistics, Johns Hopkins Bloomberg School
  of Public Health, Baltimore, MD, USA.}
\affil[3]{\small Division of Brain Injury Outcomes, Johns Hopkins Medical Institutions, Baltimore, MD, USA}
\fi

\if0\blind
\author[1]{\vspace{-2cm} }
\fi

\begin{document}\maketitle







\begin{abstract}
  We present a new estimator of the restricted mean survival time
  in randomized trials where there is
 right censoring that may depend on  treatment and baseline variables. The proposed
  estimator
  leverages prognostic baseline variables to obtain
  equal or better asymptotic precision compared to traditional estimators. Under regularity conditions and random censoring within strata of treatment and baseline variables, the proposed estimator
  has the following features: (i) it is
  interpretable under violations of the proportional hazards
  assumption; (ii) it is consistent and at least as precise as the
  Kaplan-Meier estimator under independent censoring;
  (iii) it remains consistent under violations of independent censoring (unlike the Kaplan-Meier estimator)  when either the censoring or survival distributions are estimated consistently; 
  and (iv) it achieves the nonparametric
  efficiency bound when both of these distributions are consistently estimated. We illustrate the performance of our method
  using simulations based on resampling data from a
  completed, phase 3 randomized clinical trial of a new surgical treatment for stroke; the proposed estimator achieves a 12\% gain in relative efficiency compared
  to the Kaplan-Meier estimator. The proposed estimator
  has potential advantages over existing approaches for randomized trials with time-to-event outcomes, since existing methods
  either rely on model assumptions that are untenable in many applications, or lack some of the efficiency and consistency properties (i)-(iv). We focus on estimation of the
  restricted mean survival time, but our methods may be adapted to estimate any
  treatment effect measure defined as a smooth contrast between the survival
  curves for each study arm. We provide R code to implement the
  estimator. 
\end{abstract}

%
%

\keywords{Covariate adjustment; Efficiency; Targeted minimum loss based estimation; Random censoring.}

\section{Introduction}
A standard approach to analyze clinical trials with survival outcomes
is to estimate the survival curve in each study arm using the
Kaplan-Meier estimator. This approach assumes  that the censoring time is independent of the event time for each study arm. This assumption would typically be false, e.g., if high baseline disease severity is prognostic for both earlier drop-out and earlier time to death. To accommodate this type of situation, we focus throughout on the
weaker assumption of random censoring, defined as censoring
being independent of the event time within strata of the study arm
 and baseline variables. This assumption allows informative censoring, i.e., censoring correlated with the event time in a  manner fully explained by study arm and baseline variables.

The Kaplan-Meier estimator is an unadjusted estimator, that is, it ignores baseline
variables.
Unadjusted estimators have the following potential
drawbacks: (a) they can yield inconsistent estimators of the survival
function under informative censoring as discussed above, and (b) even under independent censoring (defined as censoring independent of event time and baseline variables, for each arm), they are inefficient if
 baseline variables are prognostic of the outcome. Appropriate
adjustment for baseline variables provides an opportunity to avoid
these drawbacks by (a) providing consistent estimators of the treatment
effect under random censoring and consistent estimation of either the censoring or survival distribution, and (b) increasing efficiency (i.e., asymptotic precision) of the estimators,
thereby decreasing the required sample size and saving resources.

A commonly used method to analyze randomized trials with survival outcomes is the
proportional hazards model \citep{Cox1972}. 
 A drawback of this approach is that the treatment effect
estimate becomes uninterpretable and may be misleading under
violations of the proportional hazards assumption \citep{schemper1992cox,Tian01042014}. These violations are not necessarily easy
to detect, and can lead to false conclusions in an otherwise
well designed and executed randomized trial.

Several alternatives are available to define the effect of assignment to treatment versus control
on a survival outcome.  We focus on estimation of
the marginal (i.e., unconditional) treatment effect defined as the difference between the restricted mean
survival time (RMST) in the two study arms.  The RMST is the
expected survival time restricted to (i.e., truncated at) a time $\tau$. This parameter has
a model-free, clinically meaningful interpretation \citep{Chen2001,
  Royston2011, Zhao2012, Tian01042014}.
   For example, if $\tau=180$ days, then a 14 day improvement in RMST due to treatment means 2 more weeks alive on average during the first 6 months; this may be more directly interpretable than a hazard ratio.

Although we focus on the RMST, our methods may be adapted to estimation of any smooth
contrast between the marginal survival curves for each study arm. For example, the difference between the restricted median survival times  may be of interest, e.g., for heavy-tailed survival times.
An alternative goal, not considered here, is to estimate a
conditional effect of treatment, i.e., a contrast between
distributions conditioned on the values of certain baseline
variables. Though our ultimate goal is to estimate unconditional treatment effects, we
harness information in baseline variables both to handle informative censoring (which could lead to bias if ignored) and to adjust for  chance imbalances between study arms (to improve precision).

Our methods assume the outcome is
observed on a discrete time scale. If the outcome is measured on a
continuous time scale, our proposal may still be used by finely discretizing time. In our motivating example, we observe time at the day level, and consider a period of $\tau=180$ days.

Existing estimators for our problem can be broken into the following  three families:
methods based on modeling the survival function conditional on treatment and
baseline variables (referred to as outcome regression models), methods based on modeling the
censoring probability conditional on treatment and baseline variables, and methods called doubly
robust estimators that
combine these models. We describe these in Section~\ref{sec:litrev} where we show existing methods lack one or more of the features (i)-(iv) in the abstract. 

We propose a new doubly robust estimator of the RMST
 with all the properties in the abstract, which is derived using the general targeted minimum loss based estimation (TMLE) framework of \cite{vanderLaan&Rubin06}. Our estimator combines key ideas from \cite{Moore2009}, \cite{Rotnitzky2012}, and \cite{Gruber2012t}, as we describe in Section~\ref{sec:beating}.
 To the best of our knowledge, our proposed estimator is the first to achieve properties (i)-(iv) simultaneously for our problem.



In Section~\ref{sec:litrev}, we
review commonly used methods for our problem and
highlight their strengths and limitations. In
Section~\ref{sec:trial}, our motivating application is presented: the
analysis of a completed Phase III randomized trial of a new treatment for stroke.  In
Sections~\ref{sec:intro} and~\ref{sec_existing_estimators}, we define our estimation problem and
present estimators from related work, respectively.
 In
Section~\ref{sec:beating}, we present our new estimator.
Simulation studies are presented in Section~\ref{sec:simula}, based on our motivating application. These simulation studies demonstrate that our estimator can lead to substantial  improvements in efficiency.
We conclude with a brief discussion and directions
of future research.

\section{Related Work}\label{sec:litrev}
Various methods have been proposed that satisfy some but not all of the properties (i)-(iv) from the abstract.
\cite{Zhang2014} propose a
method to estimate the survival function at a single time point, by using a linear working model for the outcome. A working model is defined as a statistical model used to construct an estimator, but that is not necessarily assumed to be correctly specified.
Unlike our proposal,
the method of \cite{Zhang2014} requires censoring to be independent of the event time within each study arm. In addition, unlike our estimator, their method
 requires outcome regression models to be linear (and not, e.g., logistic). \cite{Parast2014} propose an estimator that adjusts for covariates through a kernel regression of the outcome on a one-dimensional dimension reduction defined as the linear predictor of a proportional hazards model. Their estimator achieves properties (i)-(ii), but not (iii)-(iv).

\citet[Section 3]{lu2011semiparametric} give a general approach for constructing estimators using longitudinal data. Their approach, if it were applied to estimate the RMST, would have  properties (i)-(ii), but not (iii)-(iv). Unlike here, they require the censoring distribution to be known. An advantage of their estimator over our proposal is the incorporation of post-baseline covariates to adjust for  time dependent confounding (i.e., when  censoring and the event time have time-varying common causes). Though we do not address this problem, the techniques in this paper may be generalized to accommodate such scenarios.
\cite{Stitelman2011} also handle time dependent confounding for survival outcomes; their estimator has properties (i), (iii), (iv) but not (ii).

Methods based on estimating the censoring distribution
\citep{Cole2004, Xie2005, Rotnitzky2005} are consistent under correct
specification of the censoring model.  However, these
estimators are typically not as efficient as  the covariate adjusted estimators described later, and thus do not fully leverage the often expensive data collected in a clinical trial.

Under consistent estimation of the censoring
 and outcome distributions at rate faster than $n^{1/4}$, doubly robust estimators are
asymptotically efficient in the nonparametric model that only assumes treatment is assigned independent of baseline variables
 \citep{Robins92,Hubbard2000,vanderLaan2003,Moore2009,Stitelman2011}. Under
outcome regression model misspecification but correct censoring model specification, doubly robust estimators remain
consistent but they can be inefficient with variance larger
than the variance of inverse probability weighted estimators. Under
independent censoring, they can also have variance larger than
unadjusted alternatives such as the Kaplan-Meier estimator. \citep[A large sample illustration, for a related problem, is in the Web Appendix of][]{diaz2015enhanced}. This is problematic
since there is no guarantee that the effort placed in constructing  adjusted
estimators will lead to improved precision in estimation of
the treatment effect.

The efficiency theory that we develop has
roots in the work of \cite{pfanzagl1982contributions, Robins92,
  Robins&Rotnitzky&Zhao94, Bickel97, hahn1998role, Scharfstein1999R,
  Bang05}, among others who laid the foundation for locally
efficient estimation of causal effects. Their methods have been
extended to incorporate enhanced efficiency properties, e.g., by
\cite{Tan2006, vanderLaan&Rubin06, zhang2008improving,
  tsiatis2008covariate, cao2009improving, tan2010bounded,
  Rotnitzky2012, Gruber2012t}.

We use the general framework of \cite{vanderLaan&Rubin06} and
\cite{Gruber2012t} to construct a targeted minimum loss based
estimator of the RMST that satisfies the properties described in the
abstract. Targeted minimum loss based estimation (TMLE) of the effect of
treatment on binary, continuous, and time to event outcomes in
randomized trials was discussed in \cite{moore2009covariate} and
\cite{Stitelman2011}. Based on the work of \cite{Gruber2012t}, \cite{diaz2015enhanced}
proposed a TMLE for ordinal outcomes with enhanced efficiency properties analogous
to those of the estimator we propose in this manuscript.

The general estimation approach of \cite{Moore2009} has properties (i), (iii), (iv), but not (ii). Our main innovation is to enhance this approach so that it also  achieves (ii). The enhancement is not trivial to achieve, and relies on the general strategy  from \cite{Rotnitzky2012} and \cite{Gruber2012t}, whose work  builds on enhanced efficiency methods described above.

\section{Motivating Application: CLEAR III Trial}\label{sec:trial}

The \href{https://clinicaltrials.gov/ct2/show/NCT00784134}{CLEAR III
  trial} (Clot Lysis: Evaluating Accelerated Resolution of
Intraventricular Hemorrhage Phase III) is a completed Phase III
multicenter randomized controlled trial.  500
 individuals with
intraventricular hemorrhage (IVH) were  randomized
with equal allocation to receive either alteplase (treatment) or saline (control) for IVH
removal.  The primary outcome was defined as a score of $3$ or less on the modified Rankin Scale (mRS) of functional disability at 180 days (where smaller values correspond to better function). The proportions with mRS at most $3$ at 180 days were
 48\% vs. 45\%  in the
treatment vs. control arms, respectively.  A key secondary outcome was
all-cause mortality at 180 days: 18\% and 29\% of patients experienced
death by 180 days in the treatment vs. control arms, respectively.
Though covariate adaptive randomization was used to assign study arms, we ignore this in illustrating our method.

\begin{figure}[!htb]
\centering
\includegraphics[scale=0.8]{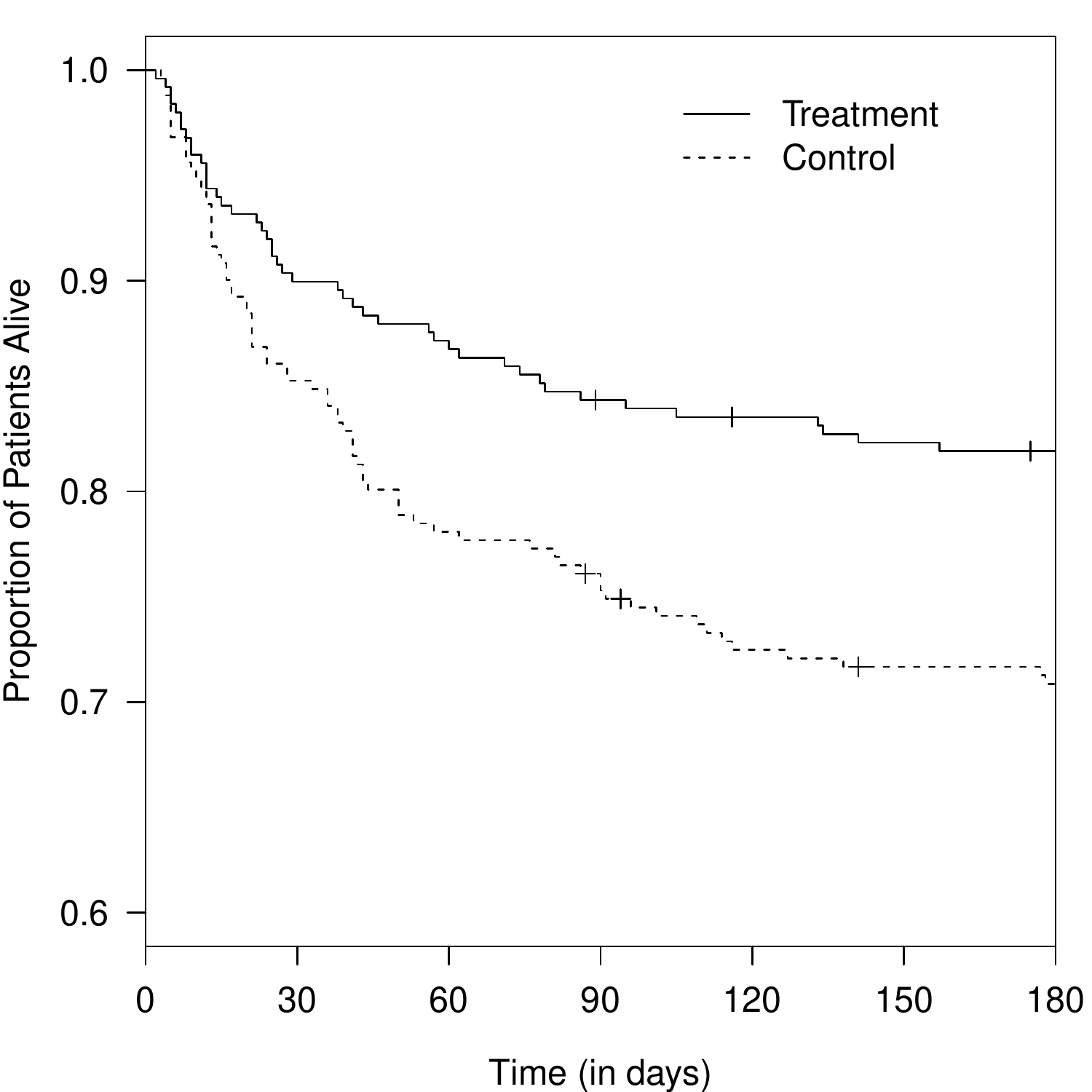}
\caption{Kaplan-Meier survival curves for the treatment arm (solid line) and control
  arm (dashed line) of the CLEAR III trial}
\label{fig:survfig}
\end{figure}

We reanalyzed data from the CLEAR III trial by defining the outcome
as time to death in days from randomization, and the treatment
effect as the difference in the RMST (at $\tau=180$ days) comparing treatment versus control arms.  Figure
~\ref{fig:survfig} displays the estimated Kaplan-Meier survival curve
for each study arm. Using the Kaplan-Meier estimator, the
unadjusted estimate of the difference in the RMST is 14.9 days
(with standard error 5.5).  The CLEAR III investigators identified several baseline
variables believed to be prognostic for   mortality; these
include age, the Glasgow Coma Score (GCS), the National Institute of
Health Stroke Scale (NIHSS) score, intracerebral hemorrhage (ICH) location (thalamus vs. other)
and ICH volume.  After adjusting for these baseline variables using
our proposed method described in Section~\ref{sec:beating}, the
estimated difference in the RMST is 14.6 days (with standard error  5.2).  Our
adjusted estimator yields an estimated variance that is roughly 12\%
smaller than the (unadjusted) RMST difference based on the
Kaplan-Meier estimator.
In the context of the CLEAR III trial, such a precision gain would allow a reduction by approximately 60 (out of original 500) patients
in the
required sample size to achieve a desired power, if a Wald-test were used based on the adjusted versus the unadjusted estimator.

\section{Data Structure, RMST Parameter, and Identification}\label{sec:intro}
\subsection{Observed Data Structure for Each Participant}
\label{subsec:observed_data} Assume $K$ equally spaced time points
$t=\{1,\dots,K\}$, e.g., representing days, at which participants are
monitored. Let $T$ denote a discrete, time-to-event outcome taking
values in $\{1,\dots,K\}\cup \{\infty\}$, where $T=\infty$ represents
no event occurring during times $1,\dots,K$.
Let $C \in \{0,\dots,K\}$ denote the censoring time defined as the
time at which the participant is last observed in the study; if a
participant remains on study through time point $K$, we let $C=K$,
which represents administrative censoring.  Let $A\in\{0,1\}$ denote
study arm assignment, and let $W$ denote a vector of baseline
variables. The observed data vector for each participant is
$O=(W,A,\Delta, \tilde T)$, where $\tilde T=\min(C,T)$, and $\Delta =
\one\{T\leq C\}$ is the indicator that the participant's event time is
observed (uncensored). Here $\one(X)$ is the indicator variable taking
value $1$ if $X$ is true and $0$ otherwise.

We assume the observed data vector for each participant $i$, denoted
$O_i=(W_i,A_i,\Delta_i, \tilde T_i)$, is an independent, identically
distributed draw from the unknown joint distribution $P_0$ on
$(W,A,\Delta, \tilde T)$. We assume $P_0\in \mathcal M$, where
$\mathcal M$ is the nonparametric model defined as all continuous
densities on $O$ with respect to a dominating measure $\nu$ such that
$A$ is independent of $W$, which holds by randomization. Our
asymptotic results are in the limit as sample size $n$ goes to
infinity, with the number of time points $K$ being fixed.

We can equivalently encode a single participant's data vector $O$ using the following longitudinal data
structure:
\begin{equation}
  O=(W, A, R_0, L_1,  R_1, L_2\ldots, R_{K-1}, L_K),\label{O}
\end{equation}
where $R_t = \one\{\tilde T = t, \Delta=0\}$ and $L_t=
\one\{\tilde T = t, \Delta=1\}$, for $t\in\{0,\ldots,K\}$. The sequence
$R_0, L_1,  R_1, L_2\ldots, R_{K-1}, L_K$ in the above display consists of all $0$'s until
the first time that either the event is observed or censoring occurs, i.e., time $t=\tilde{T}$. In the former case $L_t=1$;  otherwise $R_t=1$. For
a random variable $X$, we denote its history through time $t$ as
$\bar X_t=(X_0,\ldots,X_t)$. For a given scalar $x$, the expression
$\bar X_t=x$ denotes element-wise equality. The corresponding vector (\ref{O}) for participant $i$ is denoted by $(W_i,A_i,R_{0,i}, L_{1,i},  R_{1,i}, L_{2,i}\ldots, R_{K-1,i},L_{K,i})$.

Define the following indicator variables for each $t \geq 1$:
$$I_t=\one\{\bar R_{t-1}=0, \bar L_{t-1}=0\}, \qquad
J_t=\one\{\bar R_{t-1}=0, \bar L_t=0\}.$$
The variable $I_t$ is the indicator based on the data through time $t-1$ that a participant is
 at risk of the event being observed at time $t$; in other words, $I_t=1$ means that all the variables $R_0,L_1,R_1,L_2...,L_{t-1},R_{t-1}$ in the data vector (\ref{O}) equal $0$, which makes it possible that $L_t=1$.
 Analogously, $J_t$ is the indicator based on the outcome data  through time $t$ and censoring data before time $t$ that a participant is
 at risk of censoring at time $t$. By convention we let $J_0=1$.

Define the hazard function for survival at time $m \in \{1, \dots, K \}$:
\begin{equation}
h(m,a,w)=P_0(L_m=1|I_m = 1, A=a, W=w),\nonumber
\end{equation}
among the population at risk at
time $m$ within strata of study arm and baseline variables. Similarly, for the censoring variable $C$, define the censoring hazard at time $m \in \{0, \dots, K \}$:
\begin{equation}
g_R(m,a,w)=P_0(R_m=1|J_m=1, A=a, W=w).\nonumber
\end{equation}
We use the notation $g_A(a,w)=P_0(A=a|W=w)$ and $g=(g_A,g_R)$.
Let $p_{W}$ denote the marginal distribution of the baseline variables $W$.
We add the subscript $0$ to $p_W,g,h$ to denote the corresponding quantities under  $P_0$.
The joint distribution $P_0$ on the observed data vector $O=(W,A,\Delta, \tilde T)$ is completely characterized by the components  $p_0,g_0,h_0$, i.e.,  $P_0=(p_0,g_0,h_0)$.

\subsection{RMST Parameter Definition in Terms of Potential Outcomes}
Define the potential outcomes $T_a:a\in\{0,1\}$ as the
event times that would have been observed had study arm assignment $A=a$ and
censoring time $C=K$ been externally set with probability
one. For a restriction time $\tau\in\{1,\ldots,K\}$ of interest, the target
estimand is the difference between the restricted mean survival time setting study arm to $a=1$ versus $a=0$:
\[\theta^c = E\{\min(T_1,\tau)-\min(T_0,\tau)\}.\]
The superscript $c$ denotes a causal parameter, that is, a parameter of the distribution of the potential outcomes $T_1$ and $T_0$. We prove in the Supplementary Materials that
$E\{\min(T_a,\tau)\}=\sum_{t=0}^{\tau - 1}S^c(t,a),$ where
$S^c(t,a)=P(T_a>t)$ is the survival probability corresponding to the potential outcome
under assignment to arm $A=a$. As a
result, $\theta^c$ may be expressed as
\begin{equation}
\theta^c = \sum_{t=1}^{\tau-1}\{S^c(t,1) - S^c(t,0)\}, \label{theta_c_def}
\end{equation}
since $S^c(0,a)=1$ for $a\in \{0,1\}$.

\subsection{Identification of RMST Parameter $\theta^c$ in Terms of Observed Data Generating Distribution $P_0$} \label{sec:identification}
We show how the RMST parameter $\theta^c$, which is defined above in terms of potential outcomes, can be equivalently expressed as a function $\theta$ of the observed data distribution $P_0(W,A,\Delta, \tilde T)$, under the assumptions (a)-(d) below.
 This is useful since the potential outcomes are not always observed, in contrast to the observed data vector $(W,A,\Delta, \tilde T)$ for each participant (whose distribution we can make direct statistical inferences about); we refer to  $\theta$ as a statistical parameter, which is shorthand for saying it is a mapping from the observed data distribution $P \in \mathcal{M}$ to $\mathbb{R}$.

Define the following assumptions:
\begin{enumerate}[(a)]
\item $T= \one(A=0) T_0 + \one(A=1) T_1$ (\textit{consistency});
\item $A$ is independent of $(T_a,W)$, for each
  $a\in\{0,1\}$ (\textit{randomization});
\item $C$ is independent of $T_a$ conditional on $(A,W)$, for each
  $a\in\{0,1\}$ (\textit{random censoring});
\item $g_{A,0}(a,w) > 0$ and $g_{R,0}(t,a,w) < 1$ whenever the $P_0$-density of $W$ is positive at $W=w$, for each $a\in\{0,1\}$ and $t\in\{0,\ldots,\tau-1 \}$ (\textit{positivity assumption}).
\end{enumerate}
We make assumptions (a)-(d) throughout the manuscript.
Assumption (a) connects
 the potential outcomes to the observed outcome. Assumption (b) holds by design in a randomized trial. Assumption (c), which is similar to that in \cite{rubin1987multiple}, means that censoring is random within strata of treatment and baseline variables (which we abbreviate as ``random censoring").
 Assumption (d) states that each treatment arm has a positive probability, and that every time point has a hazard of censoring smaller than one, within each baseline variable stratum $W=w$ with positive density under $P_0$.

Denote the survival function for $T$ at time $t \in \{1,\dots, \tau-1\}$
conditioned on study arm $a$ and baseline variables $w$ by
\begin{equation}
S(t,a,w)=P(T>t|A=a,W=w).\label{S_def}
\end{equation}
Similarly, define the following function of the censoring distribution:
\begin{equation}
G(t,a,w)=P(C\geq t| A=a,W=w). \label{G_def} \end{equation}
Under assumptions (a)-(d), we have $T \indep C | A,W$ and therefore
$S(t,a,w)$ and $G(t,a,w)$  have the following product formula representations:
\begin{align}
S(t,a,w)&=\prod_{m=1}^t \{1-h(m,a,w)\}; \qquad
G(t,a,w)=\prod_{m=0}^{t-1} \{1-g_R(m,a,w)\}.
\label{defS}
\end{align}

The potential outcome survival function $S^c(t,a)$ can be equivalently represented in terms of the observed data distribution as
\begin{equation}
S(t,a)=E_{p_{W}}\prod_{m=1}^t \{1-h(m,a,W)\}, \label{DefSta}
\end{equation}
for $t\in \{1,\ldots,K\}$, $a\in\{0,1\}$;
equality of $S^c(t,a)$ and the above display follows from (\ref{defS}) and
$$S^c(t,a) = P(T_a>t)= E_{p_{W}} P(T_a>t |W) =
E_{p_{W}} P(T>t|A=a,W) = E_{P_{W}}S(t,a,W),$$
where the third equality above follows from (a) and (b).

It follows from (\ref{theta_c_def}) that the causal parameter $\theta^c$ is equal to the following statistical parameter:
\begin{equation}
\theta =
\sum_{t=1}^{\tau-1} \left[ E_{P_{W}}\left\{\prod_{m=1}^t \{1-h(m,1,W)\}\right\} - E_{P_{W}}\left\{\prod_{m=1}^t \{1-h(m,0,W)\}\right\} \right].
\label{deftheta}
\end{equation}
Our goal is to estimate $\theta$ based on $n$ independent, identically distributed observations $O_i=(W_i,A_i,\Delta_i, \tilde T_i)$ drawn from $P_0$. We  construct an estimator with properties (i)-(iv) in the abstract.

Since the parameter of interest $\theta$ is  defined as a function of $(p_{W},h)$ through (\ref{deftheta}), a natural estimation strategy would be to plug
estimates of $p_W$ and $h$ in these formulas. Estimators constructed
in this way are called substitution or plug-in
estimators, and have the advantage that they remain within bounds of
the parameter space; this is desirable in estimation of
probabilities and other bounded parameters such as the RMST. Our proposed estimator is a substitution estimator.

Define independent censoring to be $C \indep (T_a,W) \mid A$, for each $a\in\{0,1\}$. This is a stronger (more restrictive) assumption than random censoring.
We refer to a censoring mechanism $G(t,a,w)$ as non-informative if it does not depend on $w$, and as informative if it depends on $w$. Under independent censoring, $G(t,a,w)$ is non-informative.

\section{Several Estimators of RMST ($\theta$) from Related Work} \label{sec_existing_estimators}
\subsection{Unadjusted Estimators of $\theta$: Kaplan-Meier and Inverse Probability Weighted} \label{sec_unadjusted}
Throughout this subsection only, we additionally assume independent censoring.
This implies
$S(t,a)=\prod_{m=1}^t\{1-h(m,a)\}$
for $h(m,a)=P(L_m=1|I_m=1, A=a)$, as proved in the Supplementary Material.
Therefore, we have the following  simpler representation of $\theta$:
\begin{equation}
\theta = \sum_{t=1}^{\tau-1} \left[\prod_{m=1}^t \{1-h(m,1)\} -\prod_{m=1}^t \{1-h(m,0)\} \right]. \label{defKM_theta}
\end{equation}

The Kaplan-Meier estimator for $S(t,a)$ is defined as
\begin{equation}
  \Skmat =
  \prod_{m=1}^t\left\{1-\frac{\sum_{i=1}^n\one\{L_{m,i}=1, I_{m,i} =
      1, A_i=a\}}{\sum_{i=1}^n\one\{I_{m,i} = 1, A_i=a\}}\right\},\label{defkm}
\end{equation}
where we set the above fraction to be 0 if the denominator is 0.
The right side of (\ref{defkm}) was obtained by substituting the empirical counterpart of each $h(m,a)$ in the formula for $S(t,a)$ above.
The corresponding Kaplan-Meier estimator of
$\theta$ is defined analogously, as
\[\thetakm=\sum_{t=1}^{\tau-1}\{\Skmot - \Skmzt\}.\]
Since $\thetakm$ is a smooth function of at most $4(\tau-1)$ empirical means, the delta method \cite[Theorem 3.1][]{vanderVaart98} implies
\[\sqrt{n}(\thetakm - \theta) = \frac{1}{\sqrt{n}}\sum_{i=1}^n\Dkm(O_i) + o_P(1),\]
where, for an observation $O$,
\begin{equation}
  \Dkm(O) = -\sum_{t=1}^{\tau-1}\sum_{m=1}^t\left[\frac{(2A-1)I_m}{g_A(A)G(m,A)}
    \frac{S(t,A)}{S(m,A)}\left\{L_m - h(m,A)\right\}\right]
\label{defStilde}
\end{equation}
is the influence function of $\thetakm$. The above influence function $\Dkm$ may be derived from Lemma~\ref{lemmaeif} below, noting that the Kaplan-Meier
estimator is the maximum likelihood estimator in the nonparametric model where only $(A, \Delta ,\tilde T)$ (and not $W$) are measured, and therefore it is asymptotically linear with influence
function equal to the efficient influence function in the model for $(A, \Delta ,\tilde T)$. As a consequence, $\sqrt{n}(\thetakm - \theta)$
converges to a mean zero normal distribution with asymptotic variance $\mbox{Var}(\Dkm(O))$.

This estimator is
consistent and efficient for the case where only $(A,\Delta,\tilde T)$ is
observed. However, the unadjusted estimator is generally inefficient
for the case where  $(W,A,\Delta,\tilde T)$ is observed. Intuitively, this is
because the unadjusted estimator fails to leverage the prognostic
information in baseline variables $W$. Furthermore, under the less
restrictive random censoring assumption, the unadjusted
estimator will generally not be consistent, while adjusted estimators remain consistent if the censoring distribution is consistently estimated.

We next define the unadjusted, inverse probability weighted (IPW) estimator.
Let $\Gkmat$ denote the Kaplan-Meier estimator of the censoring
distribution $G(t,a)=\prod_{m=0}^{t-1}\{1-g_R(m,a)\}$, defined as
\begin{equation*}
  \Gkmat =
  \prod_{m=0}^{t-1}\left[1-\frac{\sum_{i=1}^n\one\{R_{m,i}=1, J_{m,i} =
      1, A_i=a\}}{\sum_{i=1}^n\one\{J_{m,i} = 1, A_i=a\}}\right].
\end{equation*}
The unadjusted IPW estimator of $S(t,a)$ is defined as
\[\SIPWat=\frac{1}{n}\sum_{i=1}^n \frac{\one\{A_i=a,\bar R_{t-1,i}=0,\bar L_{t,i}=0\}}{\hat
  g_A(a)\Gkmat},\]
where $\hat g_A(a)$ denotes the sample mean of $\one\{A=a\}$. \citet{tian2013restricted}
show that
\[\sqrt{n}\{\SIPWat-\Skmat\}=o_P(1),\]
i.e., the asymptotic distributions of these two estimators are equal up to $o_P(1/\sqrt{n})$, and the estimator of $\theta$ given by
$\thetaIPW=\sum_{t=1}^{\tau-1}\{\SIPWot-\SIPWzt\}$ also satisfies
\[\sqrt{n}(\thetaIPW - \theta) = \frac{1}{\sqrt{n}}\sum_{i=1}^n\Dkm(O_i) + o_P(1).\]
An important consequence is that both
$  \sqrt{n}(\thetakm - \theta)$ and $\sqrt{n}(\thetaIPW - \theta)$ converge in distribution to  $N(0,\sigmakm^2)$, where $\sigmakm^2 = \mbox{Var}(\Dkm(O))$.

\subsection{Covariate Adjustment: Inverse Probability Weighted
  Estimators}\label{sec:other}


Consider an estimator $\hat h(t, a, w)$ for the outcome hazard function $h(t, a, w)$, for example, based on fitting a Cox proportional hazards model that conditions on $A=a$ and $W=w$, and that uses the Nelson-Aalen estimator for the baseline hazard. If the proportional hazards model
 is correct, the corresponding  substitution estimator based on (\ref{deftheta}) is consistent for
$\theta$. However, under model misspecification, this substitution
estimator will generally be inconsistent. This is particularly
problematic for randomized trials with independent censoring, in
which an unadjusted, consistent estimator can be obtained through the
Kaplan-Meier survival function.

As an alternative, an adjusted, inverse probability weighted (IPW) estimator of $S(t,a)$ (where here and below we use the definition of $S(t,a)$ in (\ref{DefSta}))
is given by
\[\SADJIPWat=\frac{1}{n}\sum_{i=1}^n \frac{\one\{A_i=a,\bar R_{t-1,i}=0,\bar
  L_{t,i}=0\}}{\hat g_A(a,W_i)\hat G(t,a,W_i)},\] where $\hat g_A(a,w)$ and $\hat
G(t,a,w)$ are estimators of $g_A(a,w)$ and $G(t,a,w)$,
respectively. The adjusted IPW estimator of
the $\tau$-restricted mean survival time, given by
$\thetaADJIPW=\sum_{t=1}^{\tau-1}\{\SADJIPWot-\SADJIPWzt\}$, is consistent for
$\theta$ if the estimators $(\hat g_A, \hat g_R, \hat h)$ are consistent
at a fast enough rate. 

We focus on the case where $\hat g_R(t,a,w)$ and $\hat g(a,w)$ are
estimators having \textit{saturated terms} for time, treatment, and
their interaction. We define such estimators as estimators that
satisfy
\begin{equation*}
  \scalebox{0.9}{$\sum_{i=1}^nJ_{m,i}\{R_{m,i} -\hat g_R(m, A_i,
  W_i)\}=\sum_{i=1}^nA_iJ_{m,i}\{R_{m,i} -\hat g_R(m, A_i,
  W_i)\}=\sum_{i=1}^n \{A_i-\hat g_A(1, W_i)\}=0,$}
\end{equation*}
for $m=0,\ldots, K-1$. An example of estimators with saturated terms for time, treatment and
their interaction is given by maximum likelihood estimators for logistic regression models
of $g_R(t,a,w)$ and $g_A(a,w)$, respectively, that include at least an
intercept, main terms for $A$ and each time $t \in \{1,\dots,K\}$, and
interaction terms for $A$ by each time $t \in \{1,\dots,K\}$;
arbitrary and data-adaptive additional terms involving $t,a,w$ can be included in these
models.

When using $\hat G(t,a,w)$ and $\hat g(a,w)$ estimated with such
saturated terms, under independent censoring, $\SADJIPWat$ is
consistent for $S(t,a)$. If these models contain no terms involving
$w$, then $\SADJIPWat$ and $\SIPWat$ (both using the corresponding
estimators $\hat G$ and $\hat g$) are
identical. \citet{williamson2014variance} use this observation to show
that, under independent censoring, the asymptotic variance of
$\SADJIPWat$ is smaller or equal than the asymptotic variance of
$\SIPWat$; this, together with the delta method, implies
\[\sqrt{n}(\thetaADJIPW-\theta)\to N(0,\sigmaadjipw^2),\]
where $\sigmaadjipw^2\leq \sigmakm^2$.
(Throughout, $\to$ indicates convergence in distribution as $n$ goes to infinity.)
A weakness of $\thetaADJIPW$  is that under informative censoring, it will generally be inconsistent if the model for $G$ is misspecified. This motivates considering double robust estimators described in Section~\ref{sec_AIPW}.

Adjusted estimators (such as the adjusted IPW above) often involve fitting a parametric model for $g_A$, even though $g_A$ is known by design in a randomized trial.
 Intuitively, the purpose of this
model fit is to capture chance imbalances of the baseline variables $W$
between study arms for a given data set; these imbalances can then be adjusted to improve efficiency. The general theory underlying efficiency
improvements through estimation of known nuisance parameters such as
$g_A$ is presented, e.g., by \cite{Robins&Rotnitzky&Zhao94} and
\cite{vanderLaan2003}.

\subsection{Augmented Inverse Probability Weighted (Double Robust) Estimator} \label{sec_AIPW}
We start by presenting
the efficient influence function for estimation of $\theta$ in model $\cal
M$. The following lemma may be proved by applying the delta method to the
definition of $\theta$ in (\ref{deftheta}) and using the efficient influence
function of $S(t,a)$ as presented in \cite{Moore2009}:
\begin{lemma}\label{lemmaeif}
  The efficient influence function for estimating $\theta$ in the model
$\cal M$ is
\begin{equation} D(O)= \sum_{m=1}^{\tau-1}\big[I_mZ(m,A,W)\left\{L_m - h(m,A,W)\right\} +
  S(m,1,W) - S(m,0,W)\big] - \theta,\label{defD}
\end{equation}
where $Z(m,A,W)=Z_1(m,A,W)-Z_0(m,A,W)$, and
\begin{equation}
  Z_a(m,A,W)=-\sum_{t=m}^{\tau-1}\frac{\one\{A=a\}}{g_A(a,W)G(m,a,W)}
  \frac{S(t,a,W)}{S(m,a,W)}.\label{defZ}
\end{equation}
\end{lemma}
For conciseness, we suppress the dependence of $D$ on $g$ and $h$. The function $D$ has two important properties for estimation of
$\theta$. First, it is a doubly robust estimating function, i.e., for given
estimators $\hat{h}$ and $\hat{g}$ of $h$ and $g\equivReplacedByEquality (g_A,g_R)$,
respectively,
the estimator formed by solving for $\theta$ in the
following estimating equation
 is consistent if at least one of $h$ or $g$ is estimated consistently (while the other converges to a limit that may be incorrect):
\begin{multline}
 0= \sum_{i=1}^n \bigg\{\sum_{m=1}^{\tau-1}\bigg[\sum_{t=m}^{\tau-1}
  \frac{-(2A_i-1)I_{m,i}}{\hat g_A(A_i,W_i)\hat G(m,A_i,W_i)}\frac{\hat
    S(t,A_i,W_i)}{\hat
    S(m,A_i,W_i)}
  \{L_{m,i}-\hat{h}(m,A_i, W_i)\} \\ + \hat{S}(m,1, W_i) - \hat{S}(m,0, W_i) \bigg]
  - \theta\bigg\},\label{esteq}
\end{multline}
where
$  \hat{S}(m ,a, w)  = \prod_{m'=1}^m\{1-\hat h(m',a,w)\}$ and  $\hat{G}(m,a, w)  = \prod_{m'=0}^{m-1}\{1-\hat g_R(m',a,w)\}.$
 This estimator is often referred to as the
augmented IPW estimator, and we denote it by $\thetaAIPW$. The double robustness property is desirable since it guarantees that improper adjustment for
covariates through a misspecified working model for $h$  still leads to a consistent estimator of $\theta$ in randomized trials with random censoring if
$g_R$ is consistently estimated.

Second, the efficient
influence function (\ref{defD}) characterizes the information bound for estimation of $\theta$
in the model $\cal M$ \citep{Bickel97}. Specifically, under consistent
estimation of $h$ and $g$ at a sufficiently fast rate, $\thetaAIPW$ has variance smaller or equal to that of
any regular, asymptotically linear estimator of $\theta$ in $\cal M$ (a property called local efficiency);  in this case, if independent censoring holds, then $\thetaAIPW$ has equal or smaller asymptotic variance compared to
$\thetakm$. Unfortunately, under   misspecification of the model
for $h$, the estimator of $h$ will generally be inconsistent, which could  lead to
$\thetaAIPW$ having worse asymptotic efficiency than $\thetakm$ under independent censoring. In other words, the added robustness of $\thetaAIPW$  may come at the price of lower efficiency compared to $\thetakm$.
This
motivates the question of whether this added robustness can be achieved at no cost. I.e., can we construct
a doubly robust  estimator  with the added guarantee of equal or better asymptotic precision as $\thetakm$ if  independent censoring holds? We construct such an estimator below.


\section{Proposed Estimator}\label{sec:beating}


We develop a doubly robust estimator, denoted $\thetaTMLE$ that has properties (i)-(iv) of the abstract.
We will sometimes use the following modified
representation of the data set:
\begin{equation}
  \{(m,W_i,A_i,J_{m,i}, R_{m,i},I_{m+1,i},L_{m+1,i}): m = 0,\ldots,K-1; i
  =1,\ldots,n\}.\label{longform}
\end{equation}
This data set is referred to as the long form, and the
original data set
\begin{equation}
  \{(W_i,A_i,\Delta_i,\tilde T_i): i =1,\ldots,n\}\label{shortform}
\end{equation}
is referred to as the short form. An observation in the long form data set is a vector of the form
$(m,W_i,A_i,J_{m,i}, R_{m,i},I_{m+1,i},L_{m+1,i})$.
We define the following auxiliary variables:

\begin{align}
  M(W) =
  &\sum_{t=1}^{\tau-1}\left[\frac{S(t,1,W)}{g_A(1,W)}+\frac{S(t,0,W)}{g_A(0,W)}\right],\label{defM}\\
  H(m,A,W) = & -\sum_{t=m+1}^{\tau-1}\frac{(2A-1)}{g_A(A,W)}\frac{S(t,A,W)}{S(m,A,W)}\frac{1}{G(m+1,A,W)}.\label{defH}
\end{align}
These auxiliary covariates are constructed from a decomposition of the efficient influence function targeting improved efficiency of the estimators as described after Theorem~\ref{impEffcoro} below.

Our proposed estimator $\thetaTMLE$ is defined as the output of the following algorithm (which follows the TMLE template):

\begin{enumerate}[{Step }1.]
\item \textit{Initial estimators.} Obtain initial estimators  $\hat g_A$,
  $\hat g_R$, and $\hat h$ of  $g_A$, $g_R$, and $h$, respectively. Let $\hat{p}_W$ denote the empirical distribution of $W$.

\item \textit{Iteratively update estimates of $h$ and $g$.} Initialize $l=0$,
  and let $\hat h^l = \hat h$, $\hat g_A^l = \hat g_A$, $\hat g_R^l = \hat g_R$.
  \begin{enumerate}[(a)]
  \item {\em Update $\hat h^l$.}
  Let $\hat S^l, \hat G^l$ denote $S,G$ after substituting $\hat h^l,\hat g_R^l,\hat{p}_W$ for $h,g_R,p_W$ in (\ref{defS}) and (\ref{deftheta}).
  Augment each observation in the long form data
    set (\ref{longform}) by two covariates $\hat Z^l_{a}(m,A_i,W_i):a\in\{0,1\}$, where
    each $\hat Z^l_{a}(m,A_i,W_i)$ is constructed by substituting the estimates $\hat S^l, \hat G^l, \hat g_A^l$ evaluated at $(m,A_i,W_i)$ in (\ref{defZ}).  Estimate the
    parameter vector $\epsilon = (\epsilon_1, \epsilon_0)$ in the logistic hazard submodel $h^l_\epsilon$
    for $h$:
    \begin{equation}\logit h^l_\epsilon(m,a,w) = \logit \hat
      h^l(m,a,w) + \epsilon_1 \hat Z^l_1(m,a,w) + \epsilon_0
      \hat Z^l_0(m,a,w), \label{logistichazardsubmodel}
    \end{equation}
    by computing the following maximum likelihood estimator:
    {\small
      \begin{equation} \hat \epsilon = \arg\max_\epsilon
        \sum_{i=1}^n\sum_{m=1}^{\tau-1}I_{m,i}
        \log \left\{h_\epsilon^l(m,A_i,W_i)^{L_{m,i}}
          (1-h_\epsilon^l(m,A_i,W_i))^{1-L_{m,i}}\right\}.\label{likelihood}
      \end{equation}
    }

    The maximizer $\hat{\epsilon}$ can be computed using standard statistical
    software by a logistic regression of $L_{m,i}$ on the variables $Z^l_1(m,A_i,W_i)$,
    $\hat Z^l_0(m,A_i,W_i)$ among observations with $I_{m,i}=1$ and $m<\tau$ in the long form data set
    (\ref{longform}), and using $\logit \hat h^l(m,A_i,W_i)$ as an offset. Define
    $\hat h^{l+1}=h_{\hat \epsilon}^l$.
  \item {\em Update $\hat g_R^l$.}  Let $\hat H^l(m,A,W)$ denote $H(m,A,W)$ with
    $\hat S^l, \hat G^l, \hat g_A^l$ substituted for $S,G,g_A$,
    respectively, in (\ref{defH}).
    Augment each observation in the long form data
    set (\ref{longform}) by $\hat H^l(m,A_i,W_i)$.
    In the long form data set, estimate the parameter
    $\gamma$ in the following logistic regression submodel for $g_R(m,a,w)$:
    \begin{equation}
   \logit g_{R, \gamma}^l(m,a,w) = \logit \hat g_R^l(m,a,w) + \gamma   \hat H^l(m,a,w), \label{aug_model_g_R}
   \end{equation}
    by logistic regression of  $R_{m,i}$ on the single covariate $\hat H^l(m,A_i,W_i)$ and with offset $
    \logit \hat g_R^l(m, A,W)$
  among observations with $m<\tau-1$ and $J_{m,i}=1$. Denote the
    corresponding maximum likelihood estimate of $\gamma$ by $\hat{\gamma}$.

  \item {\em Update $\hat g_A^l$.}  Let $\hat M^l(W)$ denote $M(W)$ with $\hat g_A^l,\hat S^l$ substituted for the corresponding components
    in (\ref{defM}).
     Augment each observation in the short form data
    set by $\hat M^l(W_i)$.
    In the short form data set, estimate the parameter $\nu$
    in the following logistic regression submodel for $g_A(1|w)$:
     \begin{equation}\logit g_{A, \nu}^l(1|w) = \logit \hat g_A^l(1|w) + \nu \hat M^l(w), \label{aug_model_g_A} \end{equation}
    by logistic regression of $A_i$ on the  covariate $\hat M^l(W_i)$
    and with offset $\logit \hat g_A^l(1|W_i)$
    among all participants $i=1,\dots,n$. Denote the corresponding
    maximum likelihood estimate of $\nu$ by $\hat{\nu}$.
  \end{enumerate}
  Define $\hat h^{l+1}=h_{\hat \epsilon}^l$,   $\hat
  g_R^{l+1}=g_{R, \hat\gamma}^l$, and $\hat g_A^{l+1}=g_{A,
    \hat\nu}^l$.
\item Update $l = l+1$ and iterate the previous step until convergence. We stop
  at the first iteration for which the sample mean of the squared differences of
  predictions based on
  $\hat h^l$, $\hat g_R^l$, $\hat g_A^l$ between step $l$ and step $l+1$ is
  smaller or equal to $10^{-4}/n$.
\end{enumerate}
Denote $\hat h^\star$, $\hat g_A^\star$, and $\hat g_R^\star$ the estimators
obtained in the last iteration of the above algorithm, and define the
 enhanced efficiency TMLE estimator of $\theta$ as
\begin{equation}
  \thetaTMLE\equivReplacedByEquality
  \sum_{t=1}^{\tau-1}   \left[\frac{1}{n}\sum_{i=1}^n\prod_{m=1}^{t}\left\{1-\hat
    h^\star(m,1,W_i)\right\} -\frac{1}{n}\sum_{i=1}^n \prod_{m=1}^{t}\left\{1-\hat h^\star(m,0,W_i)\right\}\right].
\label{defTheta}
\end{equation}
The above display is the substitution estimator of $\theta$ based on (\ref{deftheta}) where the sample means over baseline variables $W_i, i=1,\dots, n$ correspond to expectation with respect to the empirical distribution of $W$.

We next consider ways to construct the initial estimators of $h$ and $g$ for step 1 above. The outcome hazard
  function $h$ may be estimated  by
  running a prediction algorithm for the probability of  $L_m=1$ as a function
  of $A$, $W$, and $m$ among observations with $I_m=1$ in the long
  form data set (\ref{longform}). 
  The censoring hazard $g_R$ may be
  estimated by running an analogous prediction algorithm of the probability that $R_m=1$ as a function
  of $A$, $W$, and $m$ among observations with $J_m=1$.
  The treatment
  mechanism $g_A$ may be estimated by fitting a parametric model for
  the probability of $A=1$ as a function of $W$ in the short form
  data set. In a randomized trial, $g_A$ is set by design.  However,
efficiency of the TMLE can be improved by estimating $g_A$ using,
  e.g., the proportion of individuals in the treatment group, or a
  logistic regression model that contains baseline variables and an intercept
  term.

Having defined the estimation algorithm, we now present our  main results  giving conditions under
which $\thetaTMLE$ is guaranteed to be at least as efficient as $\thetaADJIPW$ and $\thetakm$.

\begin{theorem}[Equal or greater asymptotic efficiency compared to adjusted IPW estimator]\label{impEff} Assume (a)-(d), $\hat
  g_A^\star$, and $\hat g_R^\star$ are $n^{1/2}$-consistent in
  $L^2(P_0)$ norm, and $\hat h^\star$ converges to some limit $h_1$ in
  $L^2(P_0)$ norm as $n \rightarrow \infty$. Then we have the
  following convergence in distribution results: \begin{align*}
    \sqrt{n}(\thetaTMLE -\theta)&\rightarrow N(0, \sigmaadjeff^2),
    \qquad &\sqrt{n}(\thetaADJIPW -\theta)&\rightarrow N(0,
    \sigmaadjipw^2),
  \end{align*}
  where $\sigmaadjeff^2\leq \sigmaadjipw^2$. In addition, if $h_1$ equals the true $h_0$, then
  $\thetaTMLE$ achieves the semiparametric efficiency bound in $\mathcal{M}$.
\end{theorem}

The consistency rates required in our previous theorem are more
restrictive than necessary to obtain the convergence in distribution
$\sqrt{n}(\thetaTMLE -\theta)\rightarrow N(0, \sigmaadjeff^2)$. In the
Supplementary Materials we present a more general theorem, along with
its proof, which shows that this convergence holds under standard,
less restrictive, doubly robust convergence assumptions on $(\hat
g^\star, \hat h^\star)$. The following is our main result:

\begin{theorem}[Equal or greater asymptotic efficiency compared to
  Kaplan-Meier estimator, under independent
  censoring]\label{impEffcoro}
  Assume (a)-(d), independent censoring, and that $\hat h^\star$
  converges as in the first sentence of Theorem~\ref{impEff}. Assume
  also that $\hat g_A^\star$, $\hat g_R^\star$ are estimated using
  models with saturated terms for time, treatment, and their
  interaction. Then
  \begin{align*}
    \sqrt{n}(\thetaTMLE -\theta)&\rightarrow N(0, \sigmaadjeff^2), \qquad 
    &\sqrt{n}(\thetaIPW -\theta)&\rightarrow N(0, \sigmakm^2), \\
    \sqrt{n}(\thetaADJIPW -\theta)&\rightarrow N(0, \sigmaadjipw^2),
    &\sqrt{n}(\thetakm -\theta)&\rightarrow N(0, \sigmakm^2),
  \end{align*}
  where $\sigmaadjeff^2 \leq \sigmaadjipw^2 \leq \sigmakm^2$.
\end{theorem}
These results guarantee that the proposed TMLE has asymptotic variance that
never exceeds that of the Kaplan-Meier estimator, under independent censoring (when the latter is consistent). To the best of our knowledge, this is the first
estimator to achieve the properties in the previous theorems for our problem.

The algorithm above (i.e., steps 1-3 and the formula (\ref{defTheta})) that generates the TMLE $\thetaTMLE$
combines key ideas from \cite{Moore2009}, \cite{Rotnitzky2012}, and \cite{Gruber2012t}. \cite{Moore2009} present a TMLE algorithm for estimating the  survival  difference $S(t,1)-S(t,0)$ in the model $\mathcal{M}$, which involves a step similar to 2a above; this estimator has properties (i), (iii), and (iv) but not (ii). The crux of our approach to achieve property (ii) (without sacrificing the other properties) is to augment the censoring and treatment models  through steps 2b and 2c. These augmented models use the covariates  (\ref{defM})-(\ref{defH}) that were specifically constructed to achieve property (ii).

The idea of augmenting censoring and treatment models to achieve enhanced efficiency properties can be traced back at least to \cite{Robins&Rotnitzky&Zhao94}. More recently, \cite{Rotnitzky2012} built on this idea to construct an estimator that has equal or better asymptotic precision than a certain parametric family of estimators including the adjusted IPW estimator, and then \cite{Gruber2012t} showed how to do the same in the TMLE framework.

It is not trivial to determine precisely how to augment the censoring and treatment models in order to guarantee (ii) holds. We explain the intuition for how we achieved this, which uses general ideas from the above related work.
Assume the conditions in Theorem~\ref{impEff}.
First,
consider the simpler case where the censoring and treatment distributions $g_R,g_A$ are known. Define the simplified TMLE to be as
in steps 1-3 above, except omitting steps 2b and 2c and using the known $g_R,g_A$ in step 2a.
 The simplified TMLE's influence function equals the  influence function of the adjusted IPW estimator minus the following expression (derived in the Supplementary Materials):
\begin{equation}
  M(W)\left[A - g_A(1,W)\right]
  +\sum_{m=0}^{\tau-2}J_mH (m,A,W)\left[R_m - g_R(m,A,W)\right],
  \label{diff_IF}
  \end{equation}
 for $M$ and $H$ the auxiliary variables in (\ref{defM}) and (\ref{defH}), respectively.

 Second, consider the case where the censoring and treatment distributions $g_R,g_A$ are unknown and estimates $\hat{g}_R,\hat{g}_A$ are used by the TMLE that involves all of steps 1-3, called the enhanced efficiency TMLE (i.e., $\thetaTMLE$). The influence function for $\thetaTMLE$ equals that of the simplified TMLE minus the latter's projection on the tangent space $T_{g^*}$ spanned by the scores of the models used to estimate $g_R,g_A$. 
 (See \citet[Section 25.3]{vanderVaart98} for background on tangent spaces and projections as used here.)
Subtracting off such a projection is helpful since it can only decrease or leave unchanged the influence function's  variance, which equals the asymptotic variance of the estimator.
By augmenting the model for $g_R$ by $H$ as in (\ref{aug_model_g_R}) in step 2b, the corresponding score is the second term in (\ref{diff_IF}); by augmenting the model for $g_A$ by $M$ as in (\ref{aug_model_g_A}) in step 2c, the corresponding score is the first term in (\ref{diff_IF}). This implies (\ref{diff_IF}) is in the tangent space $T_{g^*}$, and therefore the influence function for $\thetaTMLE$ is orthogonal to (\ref{diff_IF}).
Combining the above argument with the last line in the previous paragraph, it follows that the influence function for $\thetaTMLE$  equals the
 adjusted IPW influence function minus the latter's
 projection  on $T_{g^*}$. Therefore, $\thetaTMLE$ has asymptotic variance at most that of the adjusted IPW estimator, which gives the main conclusion of Theorem~\ref{impEff}. Theorem~\ref{impEffcoro} then follows from the asymptotic equivalence of the adjusted IPW and Kaplan-Meier estimators under the added assumptions of independent censoring and the model for $g_R$ being saturated as described above.
  A more detailed argument that fleshes out and justifies  the above outline is in the Supplementary Materials.

In addition to the efficiency properties stated in the above results, our estimator inherits the doubly robust property of the \cite{Moore2009} TMLE (i.e., the TMLE above but without steps 2b and 2c). Under random censoring, this means that our proposal has two opportunities to achieve consistency in estimating the causal effect, in contrast to the proportional hazard model which relies exclusively on the assumption that the outcome regression is correctly specified.

Under the assumptions of Theorem~\ref{impEff}, for a consistent estimate $\hatsigmaadjeff$, a  Wald-type confidence interval $\thetaTMLE \pm z_{\alpha/2} \hatsigmaadjeff/\sqrt{n}$ is guaranteed to have $1-\alpha$ asymptotic coverage probability. If the initial estimators $\hat h$ and $\hat g_R$ in step 1 are $M$-estimators (e.g., if they are estimated through maximum likelihood in parametric working models), the nonparametric bootstrap \citep{efron1979bootstrap} may be used to obtain a consistent estimate $\hatsigmaadjeff$ \citep[see Corollary 3.1 and 3.2 in ][]{wellner1996bootstrapping}. The performance of the nonparametric bootstrap is unknown when $\hat h$ or $\hat g_R$ are data-adaptive estimators (e.g., if $\hat h$ involves variable selection). The development of a consistent variance estimator in this case remains an open question.

\section{Simulation Study}\label{sec:simula}
\subsection{Data Generating Distributions}
Our simulated distributions are based on resampling data from the CLEAR III trial (described in Section~\ref{sec:trial})
  in order to mimic key features of the trial.  The baseline variables for
each participant are $W = (W_1, W_2, W_3, W_4, W_5)$ = (age, GCS, NIHSS score,
ICH location, ICH volume), each scaled to have mean 0 and standard
deviation 1. The outcome is time until death in days, denoted by $T$.
Define the treatment effect $\theta$ to be the difference in the
RMST at restriction time $\tau=180$ days. We compare the performance of the following estimators defined above: Kaplan-Meier, unadjusted IPW, adjusted IPW, augmented IPW (AIPW), and proposed TMLE. The adjusted IPW is referred to simply as IPW.
We ran simulations at two sample sizes: $n =
500$ and $n=2000$.

We consider twelve data generating distributions, where we vary the  correlation between $T$ and $W$ (using 3 scenarios labeled A, B, C defined below),  the treatment
effect $\theta$, and
the censoring mechanism. Some of these are set to mimic features of the CLEAR III data, as described below. For each data generating distribution and sample size, we generated 10,000 simulated data sets, and report on the empirical distribution of each estimator.

The key feature that determines the magnitude of precision gains for the adjusted estimators is how correlated the baseline variables are with the outcome.
To give a rough sense of the observed correlations in the CLEAR III trial data set, we fit a
 logistic regression model for the hazard of death that includes time, treatment, time by treatment interactions,
and main effects for each baseline variable in $W = (W_1, W_2, W_3, W_4, W_5)$.
We report $\exp(\hat\beta_j):j=1,2,3,4,5$, where
$\hat\beta_j$ is the estimated coefficient for
 $W_j$ in the model fit.
These were 1.5, 0.9,
1.7, 1.0, and 1.1, respectively.
We do not assume the above  model
is correct; we only used it to roughly measure correlations in the CLEAR III data set.
In our simulations,
the data
generating distributions for scenario A mimic the above  correlations; in scenarios
B and C, the correlations are reduced as described below.

In each simulated trial, the
data vector $(W,A,T,C)$ for each participant is generated as an independent,
identically distributed draw from a joint distribution $P$ that satisfies assumptions (a)-(d). This distribution will depend on the
scenario A-C, censoring mechanism, and treatment effect $\theta$.
For scenario A, each simulated participant's data vector is   generated by first
resampling a participant  with replacement from the CLEAR III  data among the 491 patients who did not drop out, and recording  only his/her pair $(W,T)$.  The resulting distribution
preserves the empirical correlation between $W$ and $T$ from the CLEAR III trial. Such resampling may lead to more realistically  complex distributions than drawing from a regression model fit for $T$ given $W$.
In each
scenario B and C, we generate an initial $(W,T)$ as just described, and
with probability $0.5$ and $1$, respectively, we replace $W$ by an independent
draw (with replacement) from the marginal distribution of $W$ in the CLEAR III
trial. The impact is that the correlation between $T$ and $W$ decreases as one
proceeds from scenario A to scenario C, with scenario C having $T$ and $W$ independent (i.e., baseline variables not prognostic for the outcome).

Next, we assign $A$ independent of $(W,T)$ by a Bernoulli
draw with probability $0.5$ of being 1 or 0. This ensures that the randomization assumption (b) holds, and induces a distribution with $\theta=0$ (no treatment effect). We also consider  distributions with  a
positive treatment effect ($\theta = 14.9$); this value was selected since it is the unadjusted Kaplan-Meier
estimate of $\theta$ from the CLEAR III data.
To obtain such distributions, we generate  $(W,A,T)$ as above except now
if  $A = 1$ we add to $T$ an independently
generated draw from a $\chi^2$ distribution with mean $\mu = 56$, where $\mu$ was calibrated to achieve $\theta = 14.9$.

For each scenario and treatment effect $\theta \in \{0,14.9\}$, given $(W,A,T)$ we generate $C$ based on either a  non-informative or informative censoring model.
Specifically, the censoring time $C$ is generated based on drawing from the corresponding distribution $g_R(t,a,w)$ given  below, for each  $t \in \{0,\dots,T-1\}$ in turn:

\begin{itemize}
\item[] \textit{Non-Informative $g_R$: } $\logit P[R_t = 1 | J_t =
  1, A, W] = -5.5 + 0.007t$;
\item[] \textit{Informative $g_R$: } $\logit P[R_t = 1 | J_t = 1, A,
  W] = -6.5 + 0.007 t + 0.6 W_3  A + 0.3 (W_1 + W_5)$.
\end{itemize}
The percentage of censored patients due to drop-out depends on the censoring model, treatment effect and scenario.  The non-informative and informative censoring distributions yield on average 62\%--68\% and 32\%--38\%  censored patients due to drop-out, respectively (with values varying within these intervals depending on $\theta$ and the scenario A-C).   Though these are higher drop-out rates than typically expected in a randomized trial, we evaluate performance under this substantial censoring to illustrate that the estimators can have good performance even in this challenging case.
Assumptions (a)-(d) from Section~\ref{sec:identification} hold for all of our data generating distributions, while  independent censoring holds only under the non-informative censoring mechanism.

\subsection{Estimators $\hat h$, $\hat g_A$, $\hat g_R$ Used by the Adjusted Estimators}
Each of
$\hat h$, $\hat g_A$, $\hat g_R$ is based on a logistic regression working model fit.
The logistic regression model for $g_A$
includes an intercept and a main term for each component of $W$.  Since treatment
is assigned with $P(A=1|W) = 0.5$ for all simulated studies, the model for $g_A$ is correctly
specified.  To account for censoring due to patient drop-out, the model for
$g_R$ includes
saturated terms for time, treatment, and their interaction, in addition to main terms for $W_1$ and $W_5$ and a treatment by $W_3$ interaction.  Therefore, the model for $g_R$ satisfies the condition in Theorem~\ref{impEffcoro} and is correctly specified when the data is generated under either non-informative or informative censoring.
We use the long form data set to fit the model for $g_R$, and set $\hat g_R(m,a,w) = 0$ for $m$ when   the corresponding risk set is empty.
The model for $h$
consists of an intercept, main terms for time (as a real value rather than categorical) and treatment, and treatment-time interaction, as well as main terms for each component of $W$.
For scenarios A and B, the
data generating distribution is based on resampling pairs $W,T$ from the CLEAR III trial; therefore, the parametric
model for $h$ is likely to be misspecified, which could easily occur in practice.  However, in scenario C where $T$ is independent of $W$, the model
for $h$ is correctly specified.

\subsection{Simulation Results}

Tables \ref{sim:table} and \ref{sim:table2000} summarize the main results of our
simulations, at sample sizes $n=500$ and $n=2000$, respectively.
 Under non-informative censoring, where
 all the estimators are consistent, all have relatively small bias. Under
informative censoring and scenarios A and B (where the baseline covariates are prognostic for the
outcome), the Kaplan-Meier and unadjusted IPW estimators are  more biased
compared to their adjusted counterparts. The informative censoring distribution depends on baseline variables
$(W_1, W_3, W_5)$, which under scenarios A and B are also
correlated with the outcome, therefore causing confounding and non-negligible bias even at the large sample size ($n=2000$); the bias is relatively small for scenario C, since
 censoring and the outcome do not share any common causes.

We measure relative efficiency as the ratio of mean squared error (RMSE) comparing the Kaplan-Meier estimator to each of the  other estimators.
First consider scenario A, where the correlation between $W$ and $T$ mimics that from the CLEAR III trial data.
The proposed estimator $\thetaTMLE$ has relative efficiency gains in the range 12\%--14\% compared to the Kaplan-Meier estimator, under non-informative censoring.
 This means that a prespecified analysis plan using our
proposed TMLE could have required roughly $11\% \approx 1-(1/1.12)$  fewer patients (55 out of 500) to achieve the same power; alternatively, the trial sample size could be conservatively planned assuming no precision gain, and then using  $\thetaTMLE$  would lead to increased power if baseline variables are prognostic for the outcome.
The efficiency gains of  $\thetaTMLE$  are similar under both types of censoring at sample size $n=500$. However, the gains are larger under  informative censoring at $n=2000$ (reaching 53\% in one case); this is due to the bias of
 the unadjusted estimators (which, unlike the adjusted estimators, are inconsistent) making a substantial contribution to the mean squared error at the larger sample size.

Next consider scenario C, where
 baseline variables are independent of $T$. For this scenario, all of the estimators are consistent under both censoring distributions. There are small losses in relative MSE for
the adjusted estimators, due to their introducing unnecessary
 variability by adjusting for baseline variables unrelated to the outcome. These precision losses are generally smaller at the larger sample size
$n=2000$ (Table~\ref{sim:table2000}) compared to $n=500$.

The unadjusted IPW estimator has similar bias and variance as the Kaplan-Meier estimator, as predicted by theory.  The adjusted IPW estimator is not as efficient as the double robust estimators AIPW and TMLE in scenarios A and B.

The TMLE $\thetaTMLE$ has similar bias compared to the AIPW estimator. The TMLE variance is slightly smaller, which translates into a  relative MSE reduction of about 1\% comparing the TMLE to the AIPW estimator in scenarios A and B. The main advantages of the TMLE compared to the AIPW estimator are: the TMLE has property (ii), i.e., guaranteed equal or better asymptotic efficiency compared to the Kaplan-Meier estimator under independent censoring; the TMLE is guaranteed to be within the bounds of the parameter space for $\theta$. The latter property guarantees that the estimated RMST for each study arm is nonnegative, which is not guaranteed for the AIPW estimator; a negative RMST would be non-interpretable.

Using the R code in the Supplementary Materials, computation of the proposed estimator took $3.5$ and $16.7$ minutes for sample sizes $500$ and $2000$, respectively, for one of the simulated trials. This computation time includes fitting the initial estimators for the event hazard, as well as the censoring and the treatment mechanism. It was performed using R version 3.2.3 on a MacBook Air with an Intel Core i5 1.3 GHz processor and 4GB of RAM.
\begin{table}[H]
\small
  \caption{Simulation results for studies of size $n = 500$. The bias, variance (VAR), and mean squared error (MSE) are displayed for the Kaplan-Meier (KM), unadjusted inverse probability weighted (Unadj.IPW), adjusted IPW (IPW), augmented IPW (AIPW) and proposed targeted minimum loss based (TMLE) estimator.  The relative MSE (RMSE) is the ratio of  the MSE for the KM estimator to the other estimators.}
  \label{sim:table}
  \centering
  \begin{tabular}{clcccccccc}
    \hline \hline
    \multicolumn{10}{l}{Non-informative censoring} \\
                 &           & \multicolumn{4}{c}{Zero Treatment Effect} & \multicolumn{4}{c}{Positive Treatment Effect} \\
        Scenario & Estimator & Bias & Var & MSE & RMSE                   & Bias & Var & MSE & RMSE    \\
    \hline \hline
    \multirow{5}{*}{A} & KM        & 0.031  & 35.23 & 35.23 & 1.000 & 0.013  & 26.54 & 26.54 & 1.000 \\
                         & Unadj.IPW & 0.032  & 35.59 & 35.59 & 0.990 & 0.044  & 26.86 & 26.86 & 0.988 \\
                         & IPW       & 0.001  & 31.76 & 31.76 & 1.109 & 0.004  & 24.29 & 24.29 & 1.093 \\
                         & AIPW      & -0.020 & 31.22 & 31.21 & 1.128 & -0.027 & 23.95 & 23.95 & 1.108 \\
                         & TMLE      & -0.019 & 30.96 & 30.96 & 1.138 & -0.074 & 23.74 & 23.74 & 1.118 \\ \hline
    \multirow{5}{*}{B} & KM        & 0.062  & 36.09 & 36.09 & 1.000 & -0.059 & 27.40 & 27.40 & 1.000 \\
                         & Unadj.IPW & 0.063  & 36.46 & 36.46 & 0.990 & -0.029 & 27.72 & 27.72 & 0.988 \\
                         & IPW       & 0.069  & 35.89 & 35.89 & 1.005 & -0.046 & 27.26 & 27.26 & 1.005 \\
                         & AIPW      & 0.069  & 35.56 & 35.56 & 1.015 & -0.087 & 26.94 & 26.94 & 1.017 \\
                         & TMLE      & 0.069  & 35.27 & 35.27 & 1.023 & -0.133 & 26.71 & 26.73 & 1.025 \\ \hline
    \multirow{5}{*}{C} & KM        & -0.048 & 35.52 & 35.52 & 1.000 & -0.035 & 26.91 & 26.91 & 1.000 \\
                         & Unadj.IPW & -0.048 & 35.90 & 35.90 & 0.989 & -0.004 & 27.22 & 27.22 & 0.989 \\
                         & IPW       & -0.046 & 36.41 & 36.41 & 0.976 & -0.013 & 27.50 & 27.50 & 0.978 \\
                         & AIPW      & -0.046 & 36.05 & 36.05 & 0.985 & -0.048 & 27.23 & 27.23 & 0.988 \\
                         & TMLE      & -0.046 & 35.76 & 35.76 & 0.993 & -0.092 & 27.00 & 27.01 & 0.996 \\ \hline \hline
                             \multicolumn{10}{l}{Informative censoring} \\
                 &           & \multicolumn{4}{c}{Zero Treatment Effect} & \multicolumn{4}{c}{Positive Treatment Effect} \\
        Scenario & Estimator & Bias & Var & MSE & RMSE                   & Bias & Var & MSE & RMSE    \\
    \hline \hline
     \multirow{5}{*}{A} & KM        & 0.802  & 31.48 & 32.12 & 1.000 & 1.760  & 22.48 & 25.58 & 1.000 \\
                          & Unadj.IPW & 0.794  & 31.61 & 32.25 & 0.996 & 1.775  & 22.59 & 25.73 & 0.994 \\
                          & IPW       & -0.108 & 29.48 & 29.49 & 1.089 & 0.079  & 22.92 & 22.93 & 1.116 \\
                          & AIPW      & -0.045 & 29.22 & 29.22 & 1.099 & 0.097  & 22.44 & 22.44 & 1.140 \\
                          & TMLE      & -0.044 & 29.08 & 29.08 & 1.105 & 0.045  & 22.18 & 22.18 & 1.153 \\ \hline
    \multirow{5}{*}{B} & KM        & 0.368  & 33.25 & 33.38 & 1.000 & 0.882  & 23.12 & 23.90 & 1.000 \\
                         & Unadj.IPW & 0.357  & 33.34 & 33.52 & 0.996 & 0.893  & 23.23 & 24.03 & 0.994 \\
                         & IPW       & -0.084 & 33.68 & 33.68 & 0.991 & 0.046  & 24.31 & 24.31 & 0.983 \\
                         & AIPW      & -0.050 & 33.36 & 33.36 & 1.000 & 0.056  & 24.13 & 24.13 & 0.990 \\
                         & TMLE      & -0.050 & 33.09 & 33.09 & 1.009 & 0.004  & 23.93 & 23.92 & 0.999 \\ \hline
    \multirow{5}{*}{C} & KM        & 0.039  & 31.88 & 31.88 & 1.000 & 0.039  & 24.02 & 24.02 & 1.000 \\
                         & Unadj.IPW & 0.025  & 32.02 & 32.02 & 0.996 & 0.046  & 24.13 & 24.13 & 0.995 \\
                         & IPW       & 0.022  & 32.66 & 32.66 & 0.976 & 0.053  & 24.90 & 24.90 & 0.965 \\
                         & AIPW      & 0.041  & 32.51 & 32.51 & 0.980 & 0.051  & 24.77 & 24.77 & 0.970 \\
                         & TMLE      & 0.041  & 32.24 & 32.24 & 0.989 & 0.001  & 24.65 & 24.65 & 0.974 \\
    \hline \hline
  \end{tabular}
\end{table}

\begin{table}[H]
\small
  \caption{Simulation results for studies of size $n = 2000$. The bias, variance (VAR), and mean squared error (MSE) are displayed for the Kaplan-Meier (KM), unadjusted inverse probability weighted (Unadj.IPW), adjusted IPW (IPW), augmented IPW (AIPW) and proposed targeted minimum loss based (TMLE) estimator.  The relative MSE (RMSE) is the ratio of  the MSE for the KM estimator to the other estimators.}
  \label{sim:table2000}
  \centering
  \begin{tabular}{clcccccccc}
    \hline \hline
    \multicolumn{10}{l}{Non-informative censoring} \\
                 &           & \multicolumn{4}{c}{Zero Treatment Effect} & \multicolumn{4}{c}{Positive Treatment Effect} \\
        Scenario & Estimator & Bias & Var & MSE & RMSE                   & Bias & Var & MSE & RMSE    \\
    \hline \hline
\multirow{5}{*}{A} & KM        & -0.069 & 8.86 & 8.86 & 1.000 & -0.015 & 6.79 & 6.79 & 1.000 \\
                                                            & Unadj.IPW & -0.070 & 8.95 & 8.96 & 0.990 & 0.016  & 6.87 & 6.87 & 0.988 \\
                                                            & IPW       & -0.057 & 8.04 & 8.04 & 1.102 & 0.027  & 6.12 & 6.12 & 1.110 \\
                                                            & AIPW      & -0.052 & 7.88 & 7.88 & 1.125 & -0.015 & 6.00 & 6.00 & 1.133 \\
                                                            & TMLE      & -0.051 & 7.81 & 7.81 & 1.134 & -0.060 & 5.95 & 5.95 & 1.142 \\ \hline
                                      \multirow{5}{*}{B} & KM        & 0.022  & 8.92 & 8.92 & 1.000 & -0.060 & 6.68 & 6.68 & 1.000 \\
                                                            & Unadj.IPW & 0.022  & 9.01 & 9.01 & 0.990 & -0.029 & 6.76 & 6.76 & 0.989 \\
                                                            & IPW       & 0.019  & 8.77 & 8.77 & 1.018 & -0.033 & 6.60 & 6.60 & 1.012 \\
                                                            & AIPW      & 0.017  & 8.69 & 8.69 & 1.027 & -0.075 & 6.50 & 6.51 & 1.027 \\
                                                            & TMLE      & 0.017  & 8.62 & 8.62 & 1.035 & -0.119 & 6.45 & 6.46 & 1.035 \\ \hline
                                       \multirow{5}{*}{C} & KM        & -0.008 & 8.78 & 8.78 & 1.000 & -0.035 & 6.77 & 6.77 & 1.000 \\
                                                            & Unadj.IPW & -0.008 & 8.87 & 8.87 & 0.990 & -0.004 & 6.85 & 6.85 & 0.989 \\
                                                            & IPW       & -0.008 & 8.91 & 8.91 & 0.986 & -0.005 & 6.87 & 6.87 & 0.986 \\
                                                            & AIPW      & -0.009 & 8.82 & 8.82 & 0.996 & -0.048 & 6.80 & 6.80 & 0.996 \\
                                                            & TMLE      & -0.009 & 8.82 & 8.82 & 0.996 & -0.038 & 6.79 & 6.79 & 0.997 \\
                                                                \hline \hline
    \multicolumn{10}{l}{Informative censoring} \\
                 &           & \multicolumn{4}{c}{Zero Treatment Effect} & \multicolumn{4}{c}{Positive Treatment Effect} \\
        Scenario & Estimator & Bias & Var & MSE & RMSE                   & Bias & Var & MSE & RMSE    \\
    \hline \hline
  \multirow{5}{*}{A} & KM        & 0.797  & 7.78 & 8.42 & 1.000 & 1.646  & 5.65 & 8.36 & 1.000 \\
                                                           & Unadj.IPW & 0.790  & 7.82 & 8.44 & 0.997 & 1.661  & 5.68 & 8.44 & 0.991 \\
                                                            & IPW       & -0.082 & 7.32 & 7.32 & 1.150 & -0.041 & 5.66 & 5.67 & 1.478 \\
                                                            & AIPW      & -0.022 & 7.20 & 7.20 & 1.168 & -0.012 & 5.53 & 5.53 & 1.513 \\
                                                            & TMLE      & -0.023 & 7.15 & 7.15 & 1.178 & -0.057 & 5.48 & 5.48 & 1.526 \\ \hline
                                       \multirow{5}{*}{B} & KM        & 0.407  & 8.08 & 8.24 & 1.000 & 0.863  & 5.76 & 6.50 & 1.000 \\
                                                            & Unadj.IPW & 0.396  & 8.12 & 8.27 & 0.997 & 0.875  & 5.78 & 6.55 & 0.993 \\
                                                            & IPW       & -0.054 & 8.15 & 8.15 & 1.011 & 0.012  & 6.05 & 6.05 & 1.075 \\
                                                            & AIPW      & -0.016 & 8.08 & 8.08 & 1.020 & 0.020  & 5.98 & 5.98 & 1.087 \\
                                                            & TMLE      & -0.016 & 8.02 & 8.02 & 1.028 & -0.024 & 5.93 & 5.93 & 1.097 \\ \hline
                                      \multirow{5}{*}{C} & KM        & -0.012 & 8.12 & 8.12 & 1.000 & -0.006 & 5.85 & 5.85 & 1.000 \\
                                                            & Unadj.IPW & -0.025 & 8.16 & 8.16 & 0.995 & 0.002  & 5.88 & 5.88 & 0.995 \\
                                                            & IPW       & -0.028 & 8.23 & 8.23 & 0.987 & -0.010 & 5.98 & 5.98 & 0.978 \\
                                                            & AIPW      & -0.009 & 8.19 & 8.19 & 0.992 & -0.015 & 5.96 & 5.95 & 0.983 \\
                                                            & TMLE      & -0.009 & 8.12 & 8.12 & 1.000 & -0.061 & 5.90 & 5.91 & 0.991 \\
    \hline \hline
  \end{tabular}
\end{table}

\section{Discussion}\label{sec:discuss}
Under random censoring, our estimator is consistent if either the
outcome or censoring model is correctly specified, i.e., our estimator
is doubly robust. This is in contrast to the proportional hazard
model, which relies exclusively on assumptions on the outcome
model. When the dimension of the baseline variables is large relative
to the sample size, the curse of dimensionality precludes the use of
nonparametric estimators or saturated parametric models
\citep{Robins97R}. A potential way to address this is to incorporate
data-adaptive model selection in constructing the initial estimators
in step 1 of the TMLE procedure, such as model stacking
\citep{Wolpert1992} or super learning
\citep{vanderLaan&Polley&Hubbard07}. The asymptotic properties of the
resulting estimator then require conditions (f)-(g) in Theorem~3 of
the Supplementary Materials. These conditions would hold automatically
for the MLE in a parametric model, but need to be verified for
data-adaptive estimators. \cite{van2014targeted} proposed an estimator
for the case of a cross-sectional study that relaxes assumption (f),
and we conjecture that this approach might be generalizable to our
problem.

In our presentation of the TMLE algorithm we have assumed that the
initial estimators $\hat g_R$ and  $\hat g_A$ contain saturated terms
for treatment, time, and their interaction. When using data-adaptive
methods, such restriction can be avoided by including the
aforementioned saturated terms in the logistic regression models
(\ref{aug_model_g_R})-(\ref{aug_model_g_A}) in step 2 of the TMLE algorithm.

Precision gains from adjustment for prognostic baseline variables can
be converted into shorter duration trials by using information
monitoring. That is, the trial continues until a prespecified
information level is achieved. Since improved precision implies
(asymptotically) a faster information accrual rate, the trial duration
may be shorter.

Our method assumes that censoring is confounded with the time to event
only by baseline variables. In the presence of time dependent
confounding between censoring and the event time, our proposal may be
adapted by augmenting the censoring and outcome models to include
time-varying confounders.  It may be possible to retain the properties
(i)-(iv) in this context by extending the proof techniques from the
Supplementary Materials. These properties are still relevant under
time dependent confounding, since the use of proportional hazard
models often yields biased estimators in this context, as discussed by
\cite{cole2003effect}.

 We defined independent censoring as $C\indep (T_a,W)\mid A$.  Another
 possible definition of independent censoring, which is more commonly
 used when discussing the unadjusted estimator, is $C\indep T_a\mid A$
 (which leaves out $W$ altogether). The latter assumption is weaker
 (less restrictive) than the former. The covariate adjustment method
 of \cite{Zhang2014} guarantees enhanced efficiency properties under
 the latter definition. Our TMLE estimator requires the former
 definition to achieve the enhanced efficiency property (ii). On the
 other hand, our TMLE estimator has the advantage over
 \cite{Zhang2014} and $\thetakm$ of remaining consistent under
 violations to the assumption that $C\indep T_a\mid A$, as long as
 random censoring holds and at least one of the censoring or survival
 distributions is consistently estimated.

Most clinical research studies use discrete time scales to measure the
time to event. If time is measured on a continuous scale,
implementation of our methods requires discretization. The specific
choice of the discretization intervals may be guided by what is
clinically relevant. For example, in clinical applications with time
to death outcomes, the clinically relevant scale would typically be a
day. In the absence of clinical criteria to guide the choice of
discretization level, a concern is that too coarse of a discretization
may lead to potentially meaningful information losses. A question for
future research is how to optimally set the level of discretization in
order to trade off information loss versus estimator
precision. Another area for future research is to consider
discretization levels that get finer with sample
size. 

\if1\blind
\section{Acknowledgments}
This work was supported by the Patient-Centered Outcomes Research Institute [ME-1306-03198],  the U.S. Food and Drug Administration [HHSF223201400113C], and the National Institute of Neurological Disorders and Stroke (NINDS) [U01NS062851]. We thank the CLEAR III trial team.
 This work is solely the responsibility of the authors and does not represent the views of the above people and agencies.
\fi

\appendix
\section{Supplementary Materials}

\begin{itemize}
\item Proofs for the main results in the paper. (proofs.pdf file)
\item R functions to compute the proposed estimator. (functions.R file)
\end{itemize}

\bibliographystyle{plainnat}
\bibliography{TMLE}

\begin{thebibliography}{44}
\providecommand{\natexlab}[1]{#1}
\providecommand{\url}[1]{\texttt{#1}}
\expandafter\ifx\csname urlstyle\endcsname\relax
  \providecommand{\doi}[1]{doi: #1}\else
  \providecommand{\doi}{doi: \begingroup \urlstyle{rm}\Url}\fi

\bibitem[Bang and Robins(2005)]{Bang05}
Heejung Bang and James~M Robins.
\newblock Doubly robust estimation in missing data and causal inference models.
\newblock \emph{Biometrics}, 61\penalty0 (4):\penalty0 962--973, 2005.

\bibitem[Bickel et~al.(1997)Bickel, Klaassen, Ritov, and Wellner]{Bickel97}
P.J. Bickel, C.A.J. Klaassen, Y.~Ritov, and J.~Wellner.
\newblock \emph{Efficient and Adaptive Estimation for Semiparametric Models}.
\newblock Springer-Verlag, 1997.

\bibitem[Cao et~al.(2009)Cao, Tsiatis, and Davidian]{cao2009improving}
Weihua Cao, Anastasios~A Tsiatis, and Marie Davidian.
\newblock Improving efficiency and robustness of the doubly robust estimator
  for a population mean with incomplete data.
\newblock \emph{Biometrika}, 96\penalty0 (3):\penalty0 723--734, 2009.

\bibitem[Chen and Tsiatis(2001)]{Chen2001}
Pei-Yun Chen and Anastasios~A Tsiatis.
\newblock Causal inference on the difference of the restricted mean lifetime
  between two groups.
\newblock \emph{Biometrics}, 57\penalty0 (4):\penalty0 1030--1038, 2001.

\bibitem[Cole and Hern{\'a}n(2004)]{Cole2004}
Stephen~R Cole and Miguel~A Hern{\'a}n.
\newblock Adjusted survival curves with inverse probability weights.
\newblock \emph{Computer Methods and Programs in Biomedicine}, 75\penalty0
  (1):\penalty0 45--49, 2004.

\bibitem[Cole et~al.(2003)Cole, Hern{\'a}n, Robins, Anastos, Chmiel, Detels,
  Ervin, Feldman, Greenblatt, Kingsley, et~al.]{cole2003effect}
Stephen~R Cole, Miguel~A Hern{\'a}n, James~M Robins, Kathryn Anastos, Joan
  Chmiel, Roger Detels, Carolyn Ervin, Joseph Feldman, Ruth Greenblatt,
  Lawrence Kingsley, et~al.
\newblock Effect of highly active antiretroviral therapy on time to acquired
  immunodeficiency syndrome or death using marginal structural models.
\newblock \emph{American Journal of Epidemiology}, 158\penalty0 (7):\penalty0
  687--694, 2003.

\bibitem[Cox(1972)]{Cox1972}
D.~R. Cox.
\newblock Regression models and life-tables (with discussion).
\newblock \emph{Journal of the Royal Statistical Society. Series B},
  34\penalty0 (2):\penalty0 187--220, 1972.

\bibitem[D\'iaz et~al.(2015)D\'iaz, Colantuoni, and
  Rosenblum]{diaz2015enhanced}
Iv\'an D\'iaz, Elizabeth Colantuoni, and Michael Rosenblum.
\newblock Enhanced precision in the analysis of randomized trials with ordinal
  outcomes.
\newblock \emph{Biometrics}, 2015.
\newblock ISSN 1541-0420.

\bibitem[Efron et~al.(1979)]{efron1979bootstrap}
B~Efron et~al.
\newblock Bootstrap methods: Another look at the jackknife.
\newblock \emph{The Annals of Statistics}, 7\penalty0 (1):\penalty0 1--26,
  1979.

\bibitem[Gruber and van~der Laan(2012)]{Gruber2012t}
Susan Gruber and Mark~J. van~der Laan.
\newblock Targeted minimum loss based estimator that outperforms a given
  estimator.
\newblock \emph{The International Journal of Biostatistics}, 8\penalty0
  (1):\penalty0 1--22, 2012.

\bibitem[Hahn(1998)]{hahn1998role}
Jinyong Hahn.
\newblock On the role of the propensity score in efficient semiparametric
  estimation of average treatment effects.
\newblock \emph{Econometrica}, pages 315--331, 1998.

\bibitem[Hubbard et~al.(2000)Hubbard, Van Der~Laan, and Robins]{Hubbard2000}
Alan~E Hubbard, Mark~J Van Der~Laan, and James~M Robins.
\newblock Nonparametric locally efficient estimation of the treatment specific
  survival distribution with right censored data and covariates in
  observational studies.
\newblock In \emph{Statistical Models in Epidemiology, the Environment, and
  Clinical Trials}, pages 135--177. Springer, 2000.

\bibitem[Lu and Tsiatis(2011)]{lu2011semiparametric}
Xiaomin Lu and Anastasios~A Tsiatis.
\newblock Semiparametric estimation of treatment effect with time-lagged
  response in the presence of informative censoring.
\newblock \emph{Lifetime Data Analysis}, 17\penalty0 (4):\penalty0 566--593,
  2011.

\bibitem[Moore and van~der Laan(2009{\natexlab{a}})]{Moore2009}
Kelly~L. Moore and Mark~J. van~der Laan.
\newblock Increasing power in randomized trials with right censored outcomes
  through covariate adjustment.
\newblock \emph{Journal of Biopharmaceutical Statistics}, 19\penalty0
  (6):\penalty0 1099--1131, 2009{\natexlab{a}}.

\bibitem[Moore and van~der Laan(2009{\natexlab{b}})]{moore2009covariate}
Kelly~L Moore and Mark~J van~der Laan.
\newblock Covariate adjustment in randomized trials with binary outcomes:
  Targeted maximum likelihood estimation.
\newblock \emph{Statistics in Medicine}, 28\penalty0 (1):\penalty0 39--64,
  2009{\natexlab{b}}.

\bibitem[Parast et~al.(2014)Parast, Tian, and Cai]{Parast2014}
Layla Parast, Lu~Tian, and Tianxi Cai.
\newblock Landmark estimation of survival and treatment effect in a randomized
  clinical trial.
\newblock \emph{Journal of the American Statistical Association}, 109\penalty0
  (505):\penalty0 384--394, 2014.

\bibitem[Pfanzagl and Wefelmeyer(1985)]{pfanzagl1982contributions}
J~Pfanzagl and W~Wefelmeyer.
\newblock Contributions to a general asymptotic statistical theory.
\newblock \emph{Statistics \& Risk Modeling}, 3\penalty0 (3-4):\penalty0
  379--388, 1985.

\bibitem[Robins and Ritov(1997)]{Robins97R}
James~M Robins and Ya'acov Ritov.
\newblock Toward a curse of dimensionality appropriate (coda) asymptotic theory
  for semi-parametric models.
\newblock \emph{Statistics in Medicine}, 16\penalty0 (3):\penalty0 285--319,
  1997.

\bibitem[Robins and Rotnitzky(1992)]{Robins92}
J.M. Robins and A.~Rotnitzky.
\newblock Recovery of information and adjustment for dependent censoring using
  surrogate markers.
\newblock In \emph{AIDS Epidemiology}, Methodological issues. Bikh\"{a}user,
  1992.

\bibitem[Robins et~al.(1994)Robins, Rotnitzky, and
  Zhao]{Robins&Rotnitzky&Zhao94}
J.M. Robins, A.~Rotnitzky, and L.P. Zhao.
\newblock Estimation of regression coefficients when some regressors are not
  always observed.
\newblock \emph{Journal of the American Statistical Association}, 89\penalty0
  (427):\penalty0 846--866, September 1994.

\bibitem[Rotnitzky and Robins(2005)]{Rotnitzky2005}
Andrea Rotnitzky and James~M Robins.
\newblock Inverse probability weighting in survival analysis.
\newblock \emph{Encyclopedia of Biostatistics}, 2005.

\bibitem[Rotnitzky et~al.(2012)Rotnitzky, Lei, Sued, and Robins]{Rotnitzky2012}
Andrea Rotnitzky, Quanhong Lei, Mariela Sued, and James~M Robins.
\newblock Improved double-robust estimation in missing data and causal
  inference models.
\newblock \emph{Biometrika}, 99\penalty0 (2):\penalty0 439--456, 2012.

\bibitem[Royston and Parmar(2011)]{Royston2011}
Patrick Royston and Mahesh~KB Parmar.
\newblock The use of restricted mean survival time to estimate the treatment
  effect in randomized clinical trials when the proportional hazards assumption
  is in doubt.
\newblock \emph{Statistics in Medicine}, 30\penalty0 (19):\penalty0 2409--2421,
  2011.

\bibitem[Rubin(1987)]{rubin1987multiple}
Donald~B Rubin.
\newblock \emph{Multiple Imputation for Nonresponse in Surveys}.
\newblock John Wiley \& Sons, 1987.

\bibitem[Scharfstein et~al.(1999)Scharfstein, Rotnitzky, and
  Robins]{Scharfstein1999R}
Daniel~O. Scharfstein, Andrea Rotnitzky, and James~M. Robins.
\newblock Adjusting for nonignorable drop-out using semiparametric nonresponse
  models: Rejoinder.
\newblock \emph{Journal of the American Statistical Association}, 94\penalty0
  (448):\penalty0 pp. 1135--1146, 1999.
\newblock ISSN 01621459.

\bibitem[Schemper(1992)]{schemper1992cox}
Michael Schemper.
\newblock Cox analysis of survival data with non-proportional hazard functions.
\newblock \emph{Journal of the Royal Statistical Society. Series D. (The
  Statistician)}, pages 455--465, 1992.

\bibitem[Stitelman et~al.(2011)Stitelman, De~Gruttola, and van~der
  Laan]{Stitelman2011}
Ori~M Stitelman, Victor De~Gruttola, and Mark~J van~der Laan.
\newblock A general implementation of tmle for longitudinal data applied to
  causal inference in survival analysis.
\newblock \emph{The International Journal of Biostatistics}, 8\penalty0 (1),
  2011.

\bibitem[Tan(2006)]{Tan2006}
Zhiqiang Tan.
\newblock A distributional approach for causal inference using propensity
  scores.
\newblock \emph{Journal of the American Statistical Association}, 101\penalty0
  (476):\penalty0 1619--1637, 2006.

\bibitem[Tan(2010)]{tan2010bounded}
Zhiqiang Tan.
\newblock Bounded, efficient and doubly robust estimation with inverse
  weighting.
\newblock \emph{Biometrika}, 97\penalty0 (3):\penalty0 661--682, 2010.

\bibitem[Tian et~al.(2014)Tian, Zhao, and Wei]{Tian01042014}
Lu~Tian, Lihui Zhao, and L.~J. Wei.
\newblock Predicting the restricted mean event time with the subject's baseline
  covariates in survival analysis.
\newblock \emph{Biostatistics}, 15\penalty0 (2):\penalty0 222--233, 2014.
\newblock \doi{10.1093/biostatistics/kxt050}.

\bibitem[Tsiatis et~al.(2008)Tsiatis, Davidian, Zhang, and
  Lu]{tsiatis2008covariate}
Anastasios~A Tsiatis, Marie Davidian, Min Zhang, and Xiaomin Lu.
\newblock Covariate adjustment for two-sample treatment comparisons in
  randomized clinical trials: A principled yet flexible approach.
\newblock \emph{Statistics in Medicine}, 27\penalty0 (23):\penalty0 4658--4677,
  2008.

\bibitem[van~der Laan(2014)]{van2014targeted}
Mark~J van~der Laan.
\newblock Targeted estimation of nuisance parameters to obtain valid
  statistical inference.
\newblock \emph{The International Journal of Biostatistics}, 10\penalty0
  (1):\penalty0 29--57, 2014.

\bibitem[van~der Laan and Robins(2003)]{vanderLaan2003}
M.J. van~der Laan and J.M. Robins.
\newblock \emph{Unified Methods for Censored Longitudinal Data and Causality}.
\newblock Springer, New York, 2003.

\bibitem[van~der Laan and Rubin(2006)]{vanderLaan&Rubin06}
M.J. van~der Laan and D.~Rubin.
\newblock Targeted maximum likelihood learning.
\newblock \emph{The International Journal of Biostatistics}, 2\penalty0
  (1):\penalty0 Article 11, 2006.

\bibitem[van~der Laan et~al.(2007)van~der Laan, Polley, and
  Hubbard]{vanderLaan&Polley&Hubbard07}
M.J. van~der Laan, E.~Polley, and A.~Hubbard.
\newblock Super learner.
\newblock \emph{Statistical Applications in Genetics \& Molecular Biology},
  6\penalty0 (25):\penalty0 Article 25, 2007.

\bibitem[van~der Vaart(1998)]{vanderVaart98}
A.~W. van~der Vaart.
\newblock \emph{Asymptotic Statistics}.
\newblock Cambridge University Press, 1998.

\bibitem[Wellner and Zhan(1996)]{wellner1996bootstrapping}
Jon~A Wellner and Yihui Zhan.
\newblock Bootstrapping z-estimators.
\newblock \emph{University of Washington Department of Statistics Technical
  Report}, 308, 1996.

\bibitem[Williamson et~al.(2014)Williamson, Forbes, and
  White]{williamson2014variance}
Elizabeth~J Williamson, Andrew Forbes, and Ian~R White.
\newblock Variance reduction in randomised trials by inverse probability
  weighting using the propensity score.
\newblock \emph{Statistics in Medicine}, 33\penalty0 (5):\penalty0 721--737,
  2014.

\bibitem[Wolpert(1992)]{Wolpert1992}
David~H Wolpert.
\newblock Stacked generalization.
\newblock \emph{Neural Networks}, 5\penalty0 (2):\penalty0 241--259, 1992.

\bibitem[Xie and Liu(2005)]{Xie2005}
Jun Xie and Chaofeng Liu.
\newblock Adjusted kaplan--meier estimator and log-rank test with inverse
  probability of treatment weighting for survival data.
\newblock \emph{Statistics in Medicine}, 24\penalty0 (20):\penalty0 3089--3110,
  2005.

\bibitem[Zhang(2014)]{Zhang2014}
Min Zhang.
\newblock Robust methods to improve efficiency and reduce bias in estimating
  survival curves in randomized clinical trials.
\newblock \emph{Lifetime Data Analysis}, 21\penalty0 (1):\penalty0 119--137,
  2014.
\newblock \doi{10.1007/s10985-014-9291-y}.

\bibitem[Zhang et~al.(2008)Zhang, Tsiatis, and Davidian]{zhang2008improving}
Min Zhang, Anastasios~A Tsiatis, and Marie Davidian.
\newblock Improving efficiency of inferences in randomized clinical trials
  using auxiliary covariates.
\newblock \emph{Biometrics}, 64\penalty0 (3):\penalty0 707--715, 2008.

\bibitem[Zhao et~al.(2012)Zhao, Tian, Uno, Solomon, Pfeffer, Schindler, and
  Wei]{Zhao2012}
Lihui Zhao, Lu~Tian, Hajime Uno, Scott~D Solomon, Marc~A Pfeffer, Jerald~S
  Schindler, and Lee~Jen Wei.
\newblock Utilizing the integrated difference of two survival functions to
  quantify the treatment contrast for designing, monitoring, and analyzing a
  comparative clinical study.
\newblock \emph{Clinical Trials}, 9\penalty0 (5):\penalty0 570--577, 2012.

\bibitem[Zhao et~al.(2016)Zhao, Claggett, Tian, Uno, Pfeffer, Solomon, Trippa,
  and Wei]{tian2013restricted}
Lihui Zhao, Brian Claggett, Lu~Tian, Hajime Uno, Marc~A. Pfeffer, Scott~D.
  Solomon, Lorenzo Trippa, and L.~J. Wei.
\newblock On the restricted mean survival time curve in survival analysis.
\newblock \emph{Biometrics}, 72\penalty0 (1):\penalty0 215--221, 2016.
\newblock ISSN 1541-0420.
\newblock \doi{10.1111/biom.12384}.
\newblock URL \url{http://dx.doi.org/10.1111/biom.12384}.

\end{thebibliography}


\begin{thebibliography}{4}
\providecommand{\natexlab}[1]{#1}
\providecommand{\url}[1]{\texttt{#1}}
\expandafter\ifx\csname urlstyle\endcsname\relax
  \providecommand{\doi}[1]{doi: #1}\else
  \providecommand{\doi}{doi: \begingroup \urlstyle{rm}\Url}\fi

\bibitem[Robins et~al.(1994)Robins, Rotnitzky, and
  Zhao]{Robins&Rotnitzky&Zhao94}
J.M. Robins, A.~Rotnitzky, and L.P. Zhao.
\newblock Estimation of regression coefficients when some regressors are not
  always observed.
\newblock \emph{Journal of the American Statistical Association}, 89\penalty0
  (427):\penalty0 846--866, September 1994.

\bibitem[van~der Laan and Rose(2011)]{vanderLaanRose11}
M.J. van~der Laan and S.~Rose.
\newblock \emph{Targeted Learning: Causal Inference for Observational and
  Experimental Data}.
\newblock Springer, New York, 2011.

\bibitem[van~der Vaart(1998)]{vanderVaart98}
A.~W. van~der Vaart.
\newblock \emph{Asymptotic Statistics}.
\newblock Cambridge University Press, 1998.

\bibitem[Zhao et~al.(2016)Zhao, Claggett, Tian, Uno, Pfeffer, Solomon, Trippa,
  and Wei]{tian2013restricted}
Lihui Zhao, Brian Claggett, Lu~Tian, Hajime Uno, Marc~A. Pfeffer, Scott~D.
  Solomon, Lorenzo Trippa, and L.~J. Wei.
\newblock On the restricted mean survival time curve in survival analysis.
\newblock \emph{Biometrics}, 72\penalty0 (1):\penalty0 215--221, 2016.
\newblock ISSN 1541-0420.
\newblock \doi{10.1111/biom.12384}.
\newblock URL \url{http://dx.doi.org/10.1111/biom.12384}.

\end{thebibliography}

\end{document}